\documentclass[11pt]{article}

\usepackage{latexsym,amsthm,amsmath,amssymb}

\def\reals{{\mathbb R}}
\def\cplx{{\mathbb C}}
\def\N{{\mathbb N}}

\def\flf{{\mathsf{FL}_f^{(4)}}}

\def\P{{\mathbb P}}
\def\C{{\mathcal C}}
\def\O{{\mathcal O}}
\def\FL{{\mathsf{FL}}}
\def\eps{{\varepsilon}}

\newtheorem{theorem}{Theorem}[section]
\newtheorem{corollary}[theorem]{Corollary}
\newtheorem{proposition}[theorem]{Proposition}

\newtheorem{lemma}[theorem]{Lemma}

\newtheorem*{thma}{Theorem 3.11}
\newcommand{\Deg}{D} 
\makeatletter
\newcommand{\ProofEndBox}{{\ifhmode\unskip\nobreak\hfil\penalty50 \else
          \leavevmode\fi\quad\vadjust{}\nobreak\hfill$\Box$
            \finalhyphendemerits=0 \par}}
\makeatother
\newcommand{\fl}[1]{\mathsf{FL}_{#1}}
\newcommand{\pfl}[1]{\mathsf{PFL}_{#1}}
\newcommand{\proofend}{\ProofEndBox\smallskip}
\newcommand{\ignore}[1]{}
\begin{document}


\title{Incidences between points and lines in $\reals^4$}

\author{
Micha Sharir
\and
Noam Solomon
}

\maketitle

\begin{abstract}
We show that the number of incidences between $m$ distinct points
and $n$ distinct lines in $\reals^4$ is $O\left(2^{c\sqrt{\log m}}
(m^{2/5}n^{4/5}+m) + m^{1/2}n^{1/2}q^{1/4} + m^{2/3}n^{1/3}s^{1/3} +
n\right)$, for a suitable absolute constant $c$, provided that no
2-plane contains more than $s$ input lines, and no hyperplane or
quadric contains more than $q$ lines. The bound holds without the
factor $2^{c\sqrt{\log m}}$ when $m \le n^{6/7}$ or $m \ge n^{5/3}$.
Except for the factor $2^{c\sqrt{\log m}}$, the bound is tight in
the worst case.
\end{abstract}


\section{Introduction}

Let $P$ be a set of $m$ distinct points in $\reals^4$ and let $L$ be a set of
$n$ distinct lines in $\reals^4$. Let $I(P,L)$ denote the number of
incidences between the points of $P$ and the lines of $L$; that is,
the number of pairs $(p,\ell)$ with $p\in P$, $\ell\in L$, and
$p\in\ell$. If all the points of $P$ and all the lines of $L$ lie in
a common plane, then the classical Szemer\'edi--Trotter
theorem~\cite{SzT} yields the worst-case tight bound

\begin{equation} \label{inc2}
I(P,L) = O\left(m^{2/3}n^{2/3} + m + n \right) .
\end{equation}
This bound clearly also holds in $\reals^4$ (or in any other dimension),
by projecting the given lines and points onto some generic plane.
Moreover, the bound will continue to be worst-case tight by placing
all the points and lines in a common plane, in a configuration that
yields the planar lower bound.

In the recent groundbreaking paper of Guth and Katz~\cite{GK2}, an
improved bound has been derived for $I(P,L)$, for a set $P$ of $m$
points and a set $L$ of $n$ lines in $\reals^3$, provided that not
too many lines of $L$ lie in a common plane \footnote{The additional
requirement in \cite{GK2}, that no regulus contains too many lines,
is not needed for the incidence bound given below.}. Specifically, they
showed:
\begin{theorem}  [Guth and Katz~\protect{\cite{GK2}}]
\label {ttt} Let $P$ be a set of $m$ distinct points and $L$ a set
of $n$ distinct lines in $\reals^3$, and let $s\le n$ be a
parameter, such that
no plane contains more than $s$ lines of $L$. Then
\begin{equation} \label{eq:ttt}
I(P,L) = O\left(m^{1/2}n^{3/4} + m^{2/3}n^{1/3}s^{1/3} + m + n\right).
\end{equation}
This bound is tight in the worst case.
\end{theorem}

In this paper, we establish the following analogous and sharper
result in four dimensions.

\begin{theorem} \label{th:main0}
Let $P$ be a set of $m$ distinct points and $L$ a set of $n$
distinct lines in $\reals^4$, and let $q,s\le n$ be parameters,
such that (i) each hyperplane or quadric contains at most $q$ lines of $L$,
and (ii) each 2-flat contains at most $s$ lines of $L$. Then
\begin {equation}
\label {ma:in0} I(P,L) \le 2^{c\sqrt{\log m}} \left( m^{2/5}n^{4/5} + m \right)
+ A\left( m^{1/2}n^{1/2}q^{1/4} + m^{2/3}n^{1/3}s^{1/3} + n\right) ,
\end {equation}
where $A$ and $c$ are suitable absolute constants.
When $m\le n^{6/7}$ or $m\ge n^{5/3}$, we get the sharper bound
\begin {equation}
\label {ma:in0x} I(P,L) \le A \left( m^{2/5}n^{4/5} + m
+ m^{1/2}n^{1/2}q^{1/4} + m^{2/3}n^{1/3}s^{1/3} + n\right) .
\end {equation}
In general, except for the factor $2^{c\sqrt{\log m}}$, the bound is
tight in the worst case, for any values of $m,n$, and for
corresponding suitable ranges of $q$ and $s$.
\end{theorem}

The proof of Theorem~\ref{th:main0} will be by induction on $m$. To
facilitate the inductive process, we extend the theorem as follows.
We say that a hyperplane or a quadric $H$ in $\reals^4$ is
\emph{$q$-restricted} for a set of lines $L$ and for an integer
parameter $q$, if there exists a polynomial $g_H$ of degree at most
$O(\sqrt{q})$, such that each of the lines of $L$ that is contained
in $H$, except for at most $q$ lines, is contained in some
irreducible component of $H \cap Z(g_H)$ that is ruled by lines and
is not a 2-flat (see below for details). In other words, a
$q$-restricted hyperplane or quadric contains in principle at most
$q$ lines of $L$, but it can also contain an unspecified number of
additional lines, all fully contained in ruled (non-planar)
components of the zero set of some polynomial of degree
$O(\sqrt{q})$. We then have the following more general result.

\begin{theorem} \label{th:main}
Let $P$ be a set of $m$ distinct points and $L$ a set of $n$
distinct lines in $\reals^4$, and let $q$ and $s\le n$ be parameters,
such that (i') each hyperplane or quadric is $q$-restricted,
and (ii) each 2-flat contains at most $s$ lines of $L$. Then,
\begin {equation}
\label {ma:in} I(P,L) \le 2^{c\sqrt{\log m}} \left( m^{2/5}n^{4/5} + m \right)
+ A\left( m^{1/2}n^{1/2}q^{1/4} + m^{2/3}n^{1/3}s^{1/3} + n\right) ,
\end {equation}
where the parameters $A$ and $c$ are as in Theorem~\ref{th:main0}.
As in the preceding theorem, when $m\le n^{6/7}$ or $m\ge n^{5/3}$, we get the sharper bound
\begin {equation}
\label {ma:inx} I(P,L) \le A \left( m^{2/5}n^{4/5} + m
+ m^{1/2}n^{1/2}q^{1/4} + m^{2/3}n^{1/3}s^{1/3} + n\right) .
\end {equation}
Moreover, except for the factor $2^{c\sqrt{\log m}}$, the bound is tight
in the worst case, as above.
\end{theorem}

The requirement that a hyperplane or quadric $H$ be $q$-restricted
extends (i.e., is a weaker condition than) the simpler requirement
that $H$ contains at most $q$ lines of $L$. Hence,
Theorem~\ref{th:main0} is an immediate corollary of
Theorem~\ref{th:main}.

A few remarks are in order.

\noindent
(a) Only the range $\sqrt{n} \le m \le n^2$ is of interest; outside this range,
regardless of the dimension of the ambient space, we have the well known and trivial
upper bound $I(P,L)=O(m+n)$, an immediate consequence of (\ref{inc2}).

\noindent
(b) The term $m^{1/2}n^{1/2}q^{1/4}$ comes from the bound of Guth and Katz~\cite{GK2}
in three dimensions (as in Theorem~\ref{ttt}), and is unavoidable, as it can be attained
if we densely ``pack'' points and lines into hyperplanes, in patterns that realize the
bound in three dimensions within each hyperplane; see Section~\ref{sec:low} for details.

\noindent
(c) Likewise, the term $m^{2/3}n^{1/3}s^{1/3}$ comes from the planar Szemer\'edi--Trotter bound
(\ref{inc2}), and is too unavoidable, as it can be attained if we densely pack points
and lines into 2-planes, in patterns that realize the bound in (\ref{inc2});
again, see Section~\ref{sec:low}.

\noindent
(d) Ignoring these terms, and the term $n$, which is included only to cater
for the case $m<\sqrt{n}$, the two terms $m^{2/5}n^{4/5}$ and $m$ ``compete'' for dominance;
the former dominates when $m=O(n^{4/3})$ and the latter when $m=\Omega(n^{4/3})$. Thus the
bound in (\ref{ma:in}) is qualitatively different within these two ranges.

\noindent
(e) The threshold $m=n^{4/3}$ also arises in the related
problem of \emph{joints} (points incident to at least four lines not in a common hyperplane)
in a set of $n$ lines in 4-space; see \cite{KSS,Qu}, and a remark below.

By a standard argument, the theorem implies the following corollary.
\begin{corollary} \label{co:main}
Let $L$ be a set of $n$ lines in $\reals^4$, satisfying the
assumptions (i') and (ii) in Theorem \ref{th:main}, for given parameters
$q$ and $s$. Then, for any $k = \Omega(2^{c\sqrt{\log n}})$, the
number $m_{\ge k}$ of points incident to at least $k$ lines of $L$
satisfies
$$
m_{\ge k} = O\Biggl(\frac{2^{\frac43 c\sqrt{\log n}}n^{4/3}}{k^{5/3}} +
\frac{nq^{1/2}}{k^2} + \frac{ns}{k^3} + \frac{n}{k} \Biggr).
$$
\end{corollary}

\noindent{\bf Remarks.} (i) It is instructive to compare
Corollary~\ref{co:main} with the analysis of \emph{joints} in a set
$L$ of $n$ lines. In $\reals^d$, a joint of $L$ is a point incident
to at least $d$ lines of $L$, not all in a common hyperplane. As
shown in \cite{KSS,Qu}, the maximum number of joints of such a set
is $O(n^{d/(d-1)})$, and this bound is worst-case tight. In four
dimensions, this bound is $O(n^{4/3})$, which corresponds to the
numerator of the first term of the bound in Corollary~\ref{co:main}.

\noindent
(ii) The other terms cater to configurations involving co-hyperplanar or
coplanar lines. For example, when $q=n$, the second term is
$O(n^{3/2}/k^2)$, in accordance with the bound obtained in Guth and
Katz~\cite{GK2} in three dimensions, and when $s=n$, the third and
fourth terms comprise (an equivalent formulation of) the bound (\ref{inc2})
of Szemer\'edi and Trotter~\cite{SzT} for the planar case.

\noindent (iii) A major interesting and challenging problem is to
extend the bound of Corollary~\ref{co:main} for any value of $k$. In
particular, is it true that the number of intersection points of the
lines (this is the case $k=2$) is $O\left(2^{\frac 4 3 c\sqrt {\log
n}}n^{4/3}+nq^{1/2}+ns\right)$? We conjecture that this is indeed
the case.

\noindent (iv) Another challenging problem is to improve our bound,
so as to get rid of, or at least reduce the factor $2^{c\sqrt{\log
m}}$. As stated in the theorems, this can be achieved when $m\le
n^{6/7}$ or $m\ge n^{5/3}$.

Additional remarks and open issues are given in the concluding Section~\ref{sec:conc}.

\paragraph{Background.}
Incidence problems have been a major topic in combinatorial and
computational geometry for the past thirty years, starting with
the Szemer\'edi-Trotter bound \cite{SzT} back in 1983.
Several techniques, interesting in their own right, have been
developed, or adapted, for the analysis of incidences, including
the crossing-lemma technique of Sz\'ekely~\cite{Sz}, and the use of
cuttings as a divide-and-conquer mechanism (e.g., see~\cite{CEGSW}).
Connections with range searching and related problems in computational
geometry have also been noted, and studies of the Kakeya problem
(see, e.g., \cite{T}) indicate the connection between this problem and
incidence problems. See Pach and Sharir~\cite{PS} for a
comprehensive (albeit a bit outdated) survey of the topic.

The landscape of incidence geometry has dramatically changed in the
past seven years, due to the infusion, in two groundbreaking papers
by Guth and Katz~\cite{GK,GK2} (the first of which was inspired by a
similar result of Dvir~\cite{Dv} for finite fields), of new tools
and techniques drawn from algebraic geometry. Although their two
direct goals have been to obtain a tight upper bound on the number
of joints in a set of lines in three dimensions \cite{GK}, and an
almost tight lower bound for the classical distinct distances
problem of Erd{\H o}s \cite{GK2}, the new tools have quickly been
recognized as useful for incidence bounds of various sorts. See
\cite{EKS,KMSS,KMS,SSZ,SoTa,Za1,Za2} for a sample of recent works on
incidence problems that use the new algebraic machinery.

The simplest instances of incidence problems involve points and
lines. Szemer\'edi and Trotter completely solved this special case
in the plane \cite{SzT}. Guth and Katz's second paper~\cite{GK2}
provides a worst-case tight bound in three dimensions, under the
assumption that no plane contains too many lines; see
Theorem~\ref{ttt}. Under this assumption, the bound in three
dimensions is significantly smaller than the planar bound (unless
one of $m,n$ is significantly smaller than the other), and the
intuition is that this phenomenon should also show up as we move to higher
dimensions. Unfortunately, the analysis becomes more involved in
higher dimensions, and requires the development or adaptation of
progressively more complex tools from algebraic geometry.
Most of these tools still appear to be unavailable, and their
absence leads either to interesting (new) open problems in the area, or
to the need to adapt existing machinery to fit into the new context.

The present paper is a first step in this direction, which considers
the four-dimensional case. It does indeed derive a sharper, nearly
optimal bound, assuming that the configuration of points and lines
is ``truly four-dimensional'', in the precise sense spelled out in
Theorems~\ref{th:main0} and~\ref{th:main}.

We also note that studying incidence problems in four (or higher)
dimensions has already taken place in several contemporary works,
such as in Solymosi and Tao~\cite{SoTa}, Zahl~\cite{Za2}, and Basu
and Sombra~\cite{BS14} (and in work in progress by Solymosi and de
Zeeuw). These works, though, consider incidences with
higher-dimensional varieties, and the study of incidences involving
lines, presented in this paper, is new. (There are several ongoing
studies, including a companion work joint with Sheffer, that aim to
derive weaker but more general bounds involving incidences between
points and curves in higher dimensions.) For very recent related
studies, see Dvir and Gopi~\cite{DG15} and Hablicsek and
Scherr~\cite{HS14}.

Our study of point-line incidences in four dimensions has lead
us to adapt more advanced tools in algebraic geometry, such as tools
involving surfaces that are ruled by lines or by flats, including Severi's
1901 work \cite{severi}, as well as the more recent works of
Landsberg \cite{IL,Land} on osculating lines and flats to algebraic
surfaces in higher dimensions.

In a preliminary version of this study~\cite{surf-socg}, we have obtained a
weaker and more constrained bound. A discussion of the significant differences
between this preliminary work and the present one is given in the overview
of the proof, which comes next.


\paragraph {Overview of the proof.}\footnote{
In this overview we assume some familiarity of the reader with
  the new ``polynomial method'' of Guth and Katz, and with subsequent applications thereof.
  Otherwise, the overview can be skipped on first reading.}
The analysis follows the general approach of Guth and Katz~\cite{GK2},
albeit with many significant adaptations and modifications.
We use induction on $m=|P|$, but we begin the description by
ignoring this aspect (for a while).
We apply the polynomial partitioning technique of Guth and
Katz~\cite{GK2}, with some polynomial $f\in \reals[x,y,z,w]$ of
suitable degree $\Deg$, and obtain a partition of $\reals^4$ into
$O(\Deg^4)$ cells, each containing at most $O(m/\Deg^4)$ points of $P$.

In our first phase, we use
\begin{equation} \label{large}
D=O(m^{2/5}/n^{1/5}),\quad\text{ for $m=O(n^{4/3})$},\quad\quad\text{and}\quad\quad
D=O(n/m^{1/2}),\quad\text{ for $m=\Omega(n^{4/3})$} .
\end{equation}
There are three types of incidences that may arise: an incidence between a
point in some cell of the partition and a line crossing that cell, an incidence
between a point on the zero set $Z(f)$ of $f$ and a line not fully contained in $Z(f)$,
and an incidence between a point on $Z(f)$ and a line fully contained in $Z(f)$.
The above choices of $\Deg$ make it a fairly easy task to bound the number of
incidences of the first two types, and the hard part is to estimate the number
of incidences of the third kind, as we have no
control on the number of points and lines contained in $Z(f)$---in the
worst case all the points and lines could be of this kind.

At the ``other end of the spectrum,'' choosing $\Deg$ to be a constant
(as done in our preliminary aforementioned study of this problem~\cite{surf-socg} and
in other recent studies of related problems \cite{Gu14,SSZ,SoTa})
simplifies considerably the handling of incidences on $Z(f)$, but
then the analysis of incidences within the cells of the partition
becomes more involved, as the sizes of the subproblems within each
cell are too large. In the works just cited (as well as in this paper), this is handled via
induction, but the price of a naive inductive approach is three-fold:
First, the bound becomes weaker, involving additional factors of the form $O(m^\eps)$,
for any $\eps>0$ (with a constant of proportionality that depends on $\eps$).
Second, the requirement that no hyperplane or quadric contains more than $q$
lines of $L$ has to be replaced by the much more restrictive requirement, that
no variety of degree at most $c_\eps$ contains more than $q$ input lines, where
$c_\eps$ is a (fairly large) constant that depends on $\eps$ (and becomes larger as $\eps$ gets smaller).
Finally, the sharp ``lower-dimensional'' terms, such as $m^{1/2}n^{1/2}q^{1/4}$
and $m^{2/3}n^{1/3}s^{1/3}$ in our case (recall that both are worst-case tight),
do not pass through the induction
successfully, so they have to be replaced by weaker terms; see the preliminary
version~\cite{surf-socg} for such weaker terms, and \cite{SSZ} for a
similar phenomenon in a different incidence problem in three dimensions.

We note that a recent study by Guth~\cite{Gu14} reexamines the
point-line incidence problem in $\reals^3$ and presents an
alternative and simpler analysis (than the original one in
\cite{GK2}), in which he uses a constant-degree partitioning
polynomial, and manages to handle successfully the relevant
lower-dimensional term $m^{2/3}n^{1/3}s^{1/3}$ through the
induction, but the analysis still incurs the extra $m^\eps$ factors
in the bound, and needs the restrictive assumption that no algebraic
surface of some large constant maximum degree $c_\eps$ contains too
many lines. In a companion paper~\cite{SS3d}, we provide yet another
simpler derivation (which is somewhat sharper than Guth's) of such
an incidence bound in three dimensions.

Our approach is to use two different choices of the degree of the
partitioning polynomial. We first choose the large value of $D$
specified above, and show that the bound in the right-hand side of
(\ref{ma:in}) accounts for the incidences within the partition
cells, for the incidences between points on $Z(f)$ and lines not
fully contained in $Z(f)$, and for most of the cases involving
incidences between points and lines on the zero set $Z(f)$. We are
then left with ``problematic'' subsets of points and lines on
$Z(f)$, which are difficult to analyze when the degree is large.
(Informally, this happens when the lines lie in certain ruled
two-dimensional subvarieties of $Z(f)$.) To handle them, we retain
only these subsets, discard the partitioning, and start afresh with
a new partitioning polynomial of a much smaller, albeit still
non-constant degree. As the degree is now too small, we need
induction to bound the number of incidences within the partition
cells.  A major feature that makes the induction work well is that
the first partitioning step ensures that the surviving set of lines
that is passed to the induction is such that each hyperplane or
quadric is now $O(D^2)$-restricted, with respect to the set of
surviving lines, and each 2-flat contains at most $O(D)$ lines of
that set (where $D$ is the large degree used in the first
partitioning step). As a consequence, the induction works better,
and ``retains'' the lower-dimensional terms $m^{1/2}n^{1/2}q^{1/4}$
and $m^{2/3}n^{1/3}s^{1/3}$. (In fact, it does not touch them at
all, because $q$ and $s$ are not passed to the induction step.) We
still pay a small price for this approach, involving the extra
factor $2^{c\sqrt{\log m}}$ in the ``leading terms''
$m^{2/5}n^{4/5}+m$ (but not in the ``lower-dimensional'' terms).
When $m$ is ``not too close to'' $n^{4/3}$, as specified in the
theorems, induction, and the use of a second partitioning
polynomial, are not needed, and a direct analysis yields the sharper
bound in (\ref{ma:inx}), without this extra factor.

The idea of using a ``small'' degree for the partitioning polynomial
is not new, and has been applied also in \cite{SSZ,Za2}. However,
the induction process in \cite{SSZ} results in weaker
lower-dimensional terms, which we avoid here with the use of two
different partitionings. We note that we have recently applied this
approach in the aforementioned study of point-line incidences in
three dimensions~\cite{SS3d}, with a simpler analysis (than that in
\cite{Gu14,GK2}) and an improved bound (than the one in
\cite{Gu14}).

The main (and hard) part of the analysis is still in handling
incidences within $Z(f)$ in the first partitioning step, where the
degree of $f$ is large. (Similar issues arise in the second step
too, but the bounds there are generally sharper than those obtained
in the first step, simply because the degree is smaller.) This is
done as follows. We first ignore the singular points on $Z(f)$. They
will be handled separately, as points lying on the zero sets of
polynomials of smaller degree (namely, partial derivatives of $f$).
We also assume that $f$ is irreducible, by considering each
irreducible factor of the original $f$ separately (see
Section~\ref{sec:pf1} for details). This step results in a
\emph{partition} of the points of $P$ and the lines of $L$ among
several varieties, each defined by an irreducible factor of $f$ or
of some derivative of $f$, so that it suffices to bound the number
of incidences between points and lines assigned to the \emph{same}
variety. The number of ``cross-variety'' incidences is shown to be
only $O(nD)$, a bound that we are ``happy'' to pay.

We next define (a four-dimensional variant of) the \emph{flecnode
polynomial} $g:=\fl{f}^4$ of $f$ (see Salmon~\cite{salmon} for the
more classical three-dimensional variant, which is used in Guth and
Katz~\cite{GK,GK2}), which vanishes at those points $p \in Z(f)$
that are incident to a line that \emph{osculates} to $Z(f)$ (i.e.,
agrees with $Z(f)$ near $p$) up to order four (and in particular to
lines that are fully contained in $Z(f)$); see below for precise
definitions. We show that $g=\fl{f}^4$ is a polynomial of degree
$O(\Deg)$. If $g\equiv 0$ on $Z(f)$ then $Z(f)$ is ruled by lines
\footnote{That is, every point $p\in Z(f)$ is incident to a line
that is fully contained in $Z(f)$; see Salmon~\cite{Edge, GK2, Kol,
salmon, SS3dv} for definitions.} (as follows from Landsberg's
work~\cite{Land}, which provides a generalization of the classical
Cayley--Salmon theorem~\cite{GK2,salmon}). We handle this case by
first reducing it to the case where $Z(f)$ is ``infinitely ruled''
by lines, meaning that most of its points are incident to infinitely
many lines that are contained in $Z(f)$ (otherwise, we can show,
using B\'ezout's theorem, that most points are incident to at most
$6$ lines, for a total of $O(m)$ incidences), and then by using the
aforementioned result of Severi~\cite{severi} from 1901, which shows
that in this case $Z(f)$ is ruled by 2-flats (each point on $Z(f)$
is incident to a 2-flat that is fully contained in $Z(f)$), unless
$Z(f)$ is a hyperplane or a quadric. This allows us to reduce the
problem to several planar incidence problems, which are reasonably
easier to handle.

The other case is where the common zero set $Z(f,g)$ of $f$ and $g$
is two-dimensional. In this case, we decompose $Z(f,g)$ into its
irreducible components, and show that the number of incidences
between points of $P$ and lines fully contained in irreducible
components that are not 2-flats is
\begin{equation} \label{weak}
\min \left\{ O(m\Deg^2+n\Deg),\; O(m+n\Deg^4) \right\} .
\end{equation}
Both terms are too large for the standard ``large'' values of $D$,
but they are non-trivial to establish, and are useful tools for slightly
improving the bound and simplifying the analysis considerably when
$D$ is not too large---see below. The derivation of these bounds is based
on a new study of point-line incidences within ruled two-dimensional
varieties in 3-space, provided in a companion paper~\cite{SS3dv}.

The irreducible components that are 2-flats are harder to handle,
because their number can be $O(\Deg^2)$ (as follows from the
generalized version of B\'ezout's theorem~\cite{Fu84}), a number
that turns out to be too large for the purpose of our incidence
bound, when a naive analysis (with a large value of $D$) is used, so
some care is needed in this case. The difficult step in this part is
when there are many points, each contained in at least three (and in
general many) 2-flats fully contained in $Z(f,g)$ (and thus in
$Z(f)$). Non-singular points of this kind are called \emph{linearly
flat} points of $Z(f)$, naturally generalizing Guth and Katz's
notion of linearly flat points in $\reals^3$~\cite{GK2} (see also
Kaplan et al.~\cite{EKS}). Linearly flat points are also \emph{flat}
points, i.e., points where the \emph{second fundamental form} of
$Z(f)$ vanishes (e.g., see Pressley~\cite{Pr}). Flatness of a point
$p$ can be expressed, again by a suitable generalization to four
dimensions of the techniques in \cite{EKS,GK2}, by the vanishing of
nine polynomials, each of degree $\le 3\Deg-4$, at $p$, which are
constructed from $f$ and from its first and second-order
derivatives. The problem can then be reduced to the case where all
the points and lines are flat (a line is flat, when not all of its
points are singular points of $Z(f)$, and all of its non-singular
points are flat). With a careful (and somewhat intricate) probing
into the geometric properties of flat lines, we can bound the number
of incidences with flat lines by reducing the problem into several
incidence problems in three dimensions (specifically, within
hyperplanes tangent to $Z(f)$ at the flat points), and then using an
extension of Guth and Katz's bound (\ref{eq:ttt}) for each of these
problems, where, in this application, we exploit the fact that each
hyperplane contains at most $q$ lines, to obtain a better,
$q$-dependent bound.

However, as noted, the terms $O(mD^2)$ (when $n^{6/7}\le m\le
n^{4/3}$) and $O(nD^4)$ (when $n^{4/3}\le m \le n^{5/3}$) are too
large (for the choices of our ``large'' values of $D$ in
(\ref{large})). We retain and also use them in the second
partitioning step, when the degree of the partitioning polynomial is
smaller, but finesse them, for the large $D$, by showing that, after
pruning away points and lines whose incidences can be estimated
directly (within the bound (\ref{ma:inx}), not using the weaker
bounds of (\ref{weak})), we are left with subsets for which every
hyperplane or quadric is $O(D^2)$-restricted, and each 2-flat
contains at most $O(D)$ lines. However, when $m\le n^{6/7}$ or $m\ge
n^{5/3}$, the terms $O(mD^2), O(nD^4)$ are not too large, and there
is no need for this part of the analysis, and a direct application
of the bounds in (\ref{weak}) yields the sharper bound in
(\ref{ma:inx}) and simplifies the proof considerably.

For the remaining range of $m$ and $n$, we go on to our second
partitioning step. We discard $f$ and start afresh with a new
partitioning polynomial $h$ of degree $E\ll D$. As already noted,
bounding incidences within the partition cells becomes non-trivial,
and we use induction, exploiting the fact that now the parameters
$q$ and $s$ are replaced by $O(D^2)$ and $O(D)$, respectively. On
the flip side of the coin, bounding incidences within $Z(h)$ is now
simpler, because $E$ is smaller, and we can use the bounds in
(\ref{weak}) (i.e., $O(mE^2+nE)$ or $O(m+nE^4)$) to establish the
bound in (\ref{ma:in}) for the ``problematic'' incidences.

The reason for using the weaker requirement that each hyperplane and
quadric be $q$-restricted, instead of just requiring that no
hyperplane or quadric contain more than $q$ lines of $L$, is that
we do not know how to bound the overall number of lines in a hyperplane
or quadric $H$ by $O(D^2)$, because of the potential existence of
ruled components of $Z(f,g)$ within $H$, which can accommodate any
number of lines. A major difference between this case and the analysis
of ruled components in Guth and Katz's study~\cite{GK2} is that here the
overall degree of $Z(f,g)$ is $O(D^2)$, as opposed to the degree of
$Z(f)$ being only $D$ in \cite{GK2}.
This precludes the application of the techniques of Guth and Katz to our
scenario---they would lead to bounds that are too large.

We also note that our analysis of incidences within $Z(f)$ is
actually carried out (in the projective 4-space) over the complex
field, which makes it simpler, and facilitates the application of
numerous tools from algebraic geometry that are developed in this
setting. The passage from the complex projective setup back to the
real affine one is straightforward---the former is a generalization
of the latter. The real affine setup is needed only for the
construction of a polynomial partitioning, which is meaningless over
$\cplx$. Once we are within the variety $Z(f)$, we can switch to the
complex projective setup, and reap the benefits noted above.

Note that, in spite of these improvements, Theorem~\ref{th:main}
still has the peculiar feature, which is not needed in Guth and
Katz~\cite{GK2} (for the incidence bound of Theorem~\ref{ttt}),
that also requires that every \emph{quadric} be $q$-restricted
(or, in the simpler version in Theorem~\ref{th:main0}, contains
at most $q$ lines of $L$).
\footnote{This is not quite the case: Guth and Katz also require
that no \emph{regulus} contains more than $s$ (actually, $\sqrt n$)
lines, but this is made to bound the number of points incident to
just two lines, and is not needed for the incidence bound in
Theorem~\ref{ttt}.} In a recent work in progress, Solomon and
Zhang~\cite{SZ} show that this requirement cannot be dropped, by
providing a construction of a quadric that contains many points and
lines, where the number of incidences between them is significantly
larger than the bound in (\ref{ma:in}) (where $q$ now only bounds
the number of lines in a hyperplane).

\section{Algebraic Preliminaries}
\label{se:alg} In this section we collect and adapt a large part of
the machinery from algebraic geometry that we need for our analysis.
Some supplementary machinery is developed within the analysis iself.

In what follows, to facilitate the application of standard
techniques in algebraic geometry, it will be more convenient to work
over the complex field $\cplx$, and in complex projective spaces. We
do so even though Theorem~\ref{th:main} is stated (and will be
proved) only for the real affine case. The passage between the two
scenarios, in the proof of the theorem, will be straightforward, as
discussed in the preceding overview. Concretely, the realness of the
underlying field is needed only for the partitioning step itself,
which has no (simple) parallel over $\cplx$. However, after reducing
the problem to points and lines contained in $Z(f)$, it is more
convenient to carry out the analysis over $\cplx$, to allow us to
apply the algebraic machinery that we are going to present next.

\subsection{Lines on varieties} \label{se:alglineon}

We begin with several basic notions and results
in differential and algebraic geometry that we will need (see, e.g.,
Ivey and Landsberg \cite{IL}, and Landsberg~\cite{Land} for more
details). For a vector space $V$ (over $\reals$ or $\cplx$), let
$\mathbb PV$ denote its projectivization. That is, $\mathbb P V = V
\setminus \{0\} /\sim$, where $v\sim w$ iff $w=\alpha v$ for some
non-zero constant $\alpha$.

An algebraic variety is the common zero set of a finite collection
of polynomials. We call it affine, if it is defined in the affine
space, or projective, if it is defined in the projective space, in
terms of homogeneous polynomials. For an (affine) algebraic variety
$X$, and a point $p \in X$, let $T_p X$ denote the (affine) tangent
space of $X$ at the point $p$. A point $p$ is \emph{non-singular} if
$\dim T_p X = \dim X$ (see Hartshorne~\cite[Definition I.5 and
Theorem I.5.1]{Hart83}). For a point $p\in X$, let $\Sigma_p$ denote
the set of the complex lines passing through $p$ and contained in
$X$, and let $\Xi_p$ denote the union of these lines (here $X$ is
implicit in these notations). For $p$ fixed, the lines in $\Sigma_p$
can be represented by their directions, as points in $\P T_p X$. In
Hartshorne~\cite[Ex.I.2.10]{Hart83}, $\Xi_p$ is also called the
\emph{(affine) cone} over $\Sigma_p$. Clearly, $\Xi_p\subseteq T_p
X$.

Consider the special case where $X$ is a hypersurface in $\cplx^4$,
i.e., $X=Z(f)$, for a non-linear polynomial
$f\in\cplx[x,y,z,w]$, which we assume to be irreducible, where
$$
Z(f)=\{p\in \cplx^4 \mid f(p)=0\}
$$
is the \emph{zero set} of $Z(f)$.
A line $\ell_v=\{p+tv \mid t\in \cplx\}$ passing through $p$ in
direction $v$ is said to \emph{osculate} to $Z(f)$ to order $k$ at
$p$, if the Taylor expansion of $f$ around $p$ in direction $v$
vanishes to order $k$, i.e., if $f(p)=0$,
\begin {equation}
\label {eq:osc} \qquad \nabla_v f(p) = 0, \qquad  \nabla_v^2 f(p) =
0, \qquad \ldots, \qquad \nabla_v^kf(p)=0,
\end {equation}
where $\nabla_v f$ (which for uniformity we also denote as
$\nabla_v^1 f$), $\nabla_v^2 f, \ldots, \nabla_v^k f$ are,
respectively, the first, second, and higher order derivatives of
$f$, up to order $k$, in direction $v$ (where $v$ is regarded as a
vector in projective 3-space, and the derivatives are interpreted in
a scale-invariant manner---we only care whether they vanish or not).
That is, $\nabla_v f = \nabla f \cdot v, \ \nabla^2_v f = v^T H_f
v$, where $H_f$ is the \emph{Hessian matrix} of $f$, and $\nabla^i_v
f$ is similarly defined, for $i>2$, albeit with more complicated
explicit expressions. For simplicity of notation, put
$F_i(p;v):=\nabla_v^i f(p)$, for $i\ge 1$.

In fact, one can extend the definition of osculation of lines to
arbitrary varieties in any dimension (see, e.g., Ivey and
Landsberg~\cite{IL}). For a variety $X$, a point $p\in X$, and an
integer $k\ge 1$, let $\Sigma_p^k \subset \mathbb PT_p X$ denote the
variety of the lines that pass through $p$ and osculate to $X$ to
order $k$ at $p$; as before, we represent the lines in $\Sigma_p^k$,
for $p$ fixed, by their directions, as points in the corresponding
projective space. For each $k \in \N$, there is a natural inclusion
$\Sigma_p \subseteq \Sigma_p^k$.
In analogy with the previous notation, we denote by $\Xi_p^k$ the
union of the lines that pass through $p$ with directions in
$\Sigma_p^k$. We let $F(X)$ denote the variety of lines (fully)
contained in $X$; this is known as the \emph{Fano variety} of $X$,
and it is a subvariety of the $(2d-2)$-dimensional
\emph{Grassmannian manifold} of lines in $\mathbb P^d(\cplx)$; see
Harris~\cite[Lecture 6, page 63]{Har} for details, and \cite[Example
6.19]{Har} for an illustration, and for a proof that this is indeed
a variety. We will sometimes denote $F(X)$ also as $\Sigma$ (or
$\Sigma(X)$), to conform with the notation involving osculating
lines. We also let $\Sigma^k$ denote the variety of the lines
osculating to order $k$ at some point of $X$, and can be thought of
as the union of the $\Sigma_p^k$ over $p\in X$. When representing
lines in $\Sigma$ or $\Sigma^k$ we can no longer use the local
representation by directions, and instead represent them, in the
customary manner, as points within the Grassmanian manifold. Here
too $\Sigma^k$ can be shown to be a variety (within the Grassmannian
manifold) and $F(X)\subseteq \Sigma^k$ for each $k$. We also have,
for any $p\in Z$, $\Sigma_p\subseteq F(X)$ and $\Sigma_p^k\subseteq
\Sigma^k$.

\paragraph{Genericity.} We recall that a property is said to hold
\emph{generically} (or \emph{generally}) for polynomials
$f_1,\ldots, f_n$, of some prescribed degrees, if there are nonzero
polynomials $g_1,\ldots, g_k$ in the coefficients of the $f_i$'s,
such that the property holds for all $f_1,\ldots, f_n$ for which
none of the polynomials $g_j$ is zero (see, e.g., Cox et
al.~\cite[Definition 3.6]{CLO2}). In this case we say that the
collection $f_1,\ldots,f_n$ is \emph{general} or \emph{generic},
with respect to the property in question, namely, with respect to
the vanishing of the polynomials $g_1,\ldots, g_k$ that define that
property.

\subsection{Generalized B\'ezout's theorem} \label{se:algbez}


An affine (resp. projective) variety $X \subset \cplx^d$ (resp. $X
\subset \P^d(\cplx)$) is called \emph{irreducible} if, whenever $V$
is written in the form $V=V_1 \cup V_2$, where $V_1$ and $V_2$ are
affine (resp., projective) varieties, then either $V_1=V$ or
$V_2=V$.

%
%

\begin {theorem}  [Cox et al.~\protect{\cite[Theorem 4.6.2, Theorem 8.3.6]{CLO}}]
Let $V$ be an affine (resp., projective) variety. Then $V$ can be
written as a finite union
$$
V=V_1\cup \cdots \cup V_m,
$$
where $V_i$ is an irreducible affine (resp., projective) variety,
for $i=1,\ldots,m$.
\end {theorem}
If one also requires that $V_i \not\subseteq V_j$ for $i\ne j$, then
this decomposition is unique, up to a permutation (see, e.g.,
\cite[Theorem 4.6.4, Theorem 8.3.6]{CLO}), and is called the
\emph{minimal decomposition} of $V$ into irreducible components.

We next state a generalized version of B\'ezout's theorem, as given
in Fulton~\cite{Fu84}. It will be a major technical tool in our analysis.
\begin{theorem}  [Fulton \protect{\cite[Proposition 2.3]{Fu84}}]
\label{th:agbt} Let $V_1,\ldots, V_s$ be subvarieties of $\mathbb
P^d$, and let $Z_1,\ldots, Z_r$ be the irreducible components of
$\bigcap_{i=1}^s V_i$. Then
$$\sum_{i=1}^r \deg(Z_i) \le \prod_{j=1}^s \deg(V_j).$$
\end {theorem}

A simple application of Theorem~\ref{th:agbt} yields the following useful result.
\begin {lemma}
\label {le:btl} A curve $\mathcal C \subset \P^4$ of
degree $\Deg$ can contain at most $\Deg$ lines.
\end {lemma}

\noindent {\bf Proof.}
Let $t$ denote the number of these lines, and let $\C_0\subset\C$ denote their union.
Intersect $\C_0$ with a generic hyperplane $H$. By Theorem~\ref{th:agbt},
the number of intersection points satisfies
$$
t \le \deg(\C_0)\cdot \deg(H) \le \deg(\C)\cdot 1 = D ,
$$
as asserted.
\proofend

This immediately yields the following result, derived in Guth and Katz~\cite{GK}
(see also \cite{EKS}) in a somewhat different manner.
\begin {corollary}
\label {co:btl} Let $f$ and $g$ be two trivariate polynomials without a common factor.
Then $Z(f,g) := Z(f)\cap Z(g)$ contains at most $\deg(f)\cdot\deg(g)$ lines.
\end {corollary}

\noindent {\bf Proof.}
This follows since $Z(f,g)$ is a curve of degree at most $\deg(f)\cdot\deg(g)$.
\proofend

\subsection{Generically finite morphisms and the Theorem of the Fibers} \label{se:finfib}

The following results can be found, e.g., in Harris~\cite[Chapter 11]{Har}.

For a map $\pi : X \to Y$ of projective varieties, and for $y\in Y$,
the variety $\pi^{-1}(y)$ is called the \emph{fiber} of $\pi$ over $y$.

The following result is a slight paraphrasing of
Harris~\cite[Proposition 7.16]{Har} and also appear in Sharir and
Solomon~\cite[Theorem 7]{SS3dv}

\begin {theorem} [Harris~\protect{\cite[Proposition 7.16]{Har}}]
\label {co:harr} Let $f: X \to Y$ be the map induced by the standard
projection map $\pi: \P^d \to \P^{r}$ (which retains $r$ of the
coordinates and discards the rest), where $r<d$, $X \subset \P^d$
and $Y\subset \P^{r}$ are projective varieties, $X$ is irreducible,
and $Y$ is the image of $X$. Then the general fiber\footnote{%
The meaning of this statement is that the assertion holds for the
fiber at any point outside some lower-dimensional
  exceptional subvariety.}
of the map $f$ is finite if and only if $\dim(X)=\dim(Y)$. In this
case, the number of points in a general fiber of $f$ is constant.
\end {theorem}

An important technical tool for our analysis is the following so-called
Theorem of the Fibers.
\begin {theorem} [Harris~\protect{\cite[Corollary 11.13]{Har}}]
\label {th:harr} Let $X$ be a projective variety and $\pi: X \to
\P^d$ be a polynomial map (i.e., the coordinate functions $x_0\circ
\pi,\ldots,x_d \circ \pi$ are homogeneous polynomials); let
$Y=\pi(X)$ denote its image. For any $p\in Y$, let
$\lambda(p)=\dim(\pi^{-1}(p))$. Then $\lambda(p)$ is an upper
semi-continuous function of $p$ in the Zariski topology\footnote{
  The Zariski closure of a set $Y$ is the intersection of all varieties
  $X$ that contain $Y$. $Y$ is Zariski closed if it is equal to its
  closure (and is therefore a variety), and is Zariski open if its
  complement is Zariski closed. See~\cite{Hart83} for further details.}
on $Y$; that is, for any $m$, the locus of points $p\in Y$ such that
$\lambda(p)\ge m$ is closed in $Y$. Moreover, if $X_0 \subset X$ is
any irreducible component, $Y_0=\pi(X_0)$ its image, and $\lambda_0$
the minimum value of $\lambda(p)$ on $Y_0$, then
$$\dim(X_0)=\dim(Y_0)+\lambda_0.$$
\end {theorem}

\subsection{Flecnode polynomials and ruled surfaces in four dimensions}
\label{se:algruled}

\paragraph{Ruled surfaces in three dimensions.}
We first review several basic properties of ruled two-dimensional surfaces
in $\reals^3$ or in $\cplx^3$. Most of these results are considered folklore in the
literature, although we have been unable to find concrete rigorous proofs
(in the ``modern'' jargon of algebraic geometry). For the sake of
completeness we provide such proofs in a companion paper~\cite{SS3dv}.

For a modern approach to ruled surfaces, there are many references;
see, e.g., Hartshorne~\cite[Section V.2]{Hart83}, or
Beauville~\cite[Chapter III]{Beau}. We say that a real $($resp.,
complex$)$ surface $X$ is \emph{ruled by real $($}resp.,\emph{
complex$)$ lines} if every point $p\in X$ in a Zariski-open dense
set is incident to a real (resp., complex) line that is fully
contained in $X$; see, e.g., \cite{salmon} or \cite{Edge} for
further details on ruled surfaces. This definition is slightly
weaker than the classical definition, where it is required that
\emph{every} point of $X$ be incident to a line contained in $X$
(e.g., as in~\cite{salmon}). It has been used in recent works, see,
e.g.,~\cite{GK2,Kol}. Similarly to the proof of Lemma 3.4 in Guth
and Katz~\cite{GK2}, a limiting argument implies that the two
definitions are equivalent. We spell out the details in
Lemma~\ref{le:rs} in the appendix (see also Sharir and
Solomon~\cite[Lemma 11]{SS3dv}).

We note that some care has to be exercised when dealing
with ruled surfaces, because ruledness may depend on the underlying
field. Specifically, it is possible for a surface defined by real
polynomials to be ruled by complex lines, but not by real lines. For
example, the sphere defined by $x^2+y^2+z^2-1=0$, regarded as a real
variety, is certainly not ruled by lines, but as a complex variety
it is ruled by (complex) lines. (Indeed, each point $(x_0,y_0,z_0)$
on the sphere is incident to the (complex) line $(x_0+\alpha t,
y_0+\beta t, z_0 +\gamma t)$, for $t\in \cplx$, where
$\alpha^2+\beta^2+\gamma^2=0$ and $\alpha x_0 + \beta y_0 + \gamma
z_0=0$, which is fully contained in the sphere.)

In three dimensions, a two-dimensional irreducible ruled surface can be
either \emph{singly ruled}, or \emph{doubly ruled} (notions that are
elaborated below), or a plane. As the following lemma shows, the only
doubly ruled surfaces are \emph{reguli}, where a regulus is the union
of all lines that meet three pairwise skew lines. There are only two
kinds of reguli, both of which are quadrics---hyperbolic paraboloids
and hyperboloids of one sheet; see, e.g., Fuchs and Tabachnikov~\cite{FT}
for more details.

The following (folklore) lemma provides a (somewhat stronger than usual)
characterization of doubly ruled surfaces; see \cite{SS3dv} for a proof.
\begin {lemma}
\label{doubly} Let $V$ be an irreducible ruled surface in $\reals^3$ or
in $\cplx^3$ which is not a plane, and let $\C\subset V$ be an algebraic
curve, such that every non-singular point $p\in V\setminus \C$ is incident
to exactly two lines that are fully contained in $V$. Then $V$ is a regulus.
\end {lemma}

When $V$ is an irreducible ruled surface which is neither a plane
nor a regulus, it must be \emph{singly ruled}, in the precise sense
spelled out in the following theorem (see also \cite{GK2}); again,
see \cite[Theorem 10]{SS3dv} for a proof.

\begin{theorem}
\label{singly} (a) Let $V$ be an irreducible ruled two-dimensional
surface of degree $\Deg>1$ in $\reals^3$ (or in $\cplx^3$), which is
not a regulus. Then, except for at most two exceptional lines, the
lines that are fully contained in $V$ are parametrized by an
irreducible algebraic curve $\Sigma_0$ in the Pl\"ucker space
$\P^5$, and thus yield a 1-parameter family of generator lines
$\ell(t)$, for $t\in \Sigma_0$, that depend continuously on the real
or complex parameter $t$. Moreover, if $t_1 \ne t_2$, and $\ell(t_1)
\ne \ell(t_2)$, then there exist sufficiently small and disjoint
neighborhoods $\Delta_1$ of $t_1$ and $\Delta_2$ of $t_2$, such that
all the lines $\ell(t)$, for $t\in \Delta_1\cup \Delta_2$, are
distinct.

\smallskip

\noindent (b) There exists a one-dimensional curve $\C\subset V$,
such that any point $p$ in $V\setminus\C$ is incident to exactly one
generator line of $V$.
\end{theorem}


Following this theorem, we refer to irreducible ruled surfaces that
are neither planes nor reguli as \emph{singly ruled}. A line $\ell$,
fully contained in an irreducible singly ruled surface $V$, such
that every point of $\ell$ is incident to another line fully
contained in $V$, is called an \emph{exceptional line} of $V$ (these
are the lines mentioned in Theorem~\ref{singly}(a)). If there exists
a point $p_V \in V$, which is incident to infinitely many lines
fully contained in $V$, then $p_V$ is called an \emph{exceptional
point} of $V$. By Guth and Katz~\cite{GK2}, $V$ can contain at most
one exceptional point $p_V$ (in which case $V$ is a cone with $p_V$
as its apex), and (as also asserted in the theorem) at most two
exceptional lines.

\paragraph{The flecnode polynomial in four dimensions.}
Let $f \in \cplx[x,y,z,w]$ be a polynomial of degree $\Deg\ge 4$. A
\emph{flecnode} of $f$ is a point $p\in Z(f)$ for which there exists
a line that passes through $p$ and osculates to $Z(f)$ to order four
at $p$. Therefore, if the direction of the line is
$v=(v_0,v_1,v_2,v_3)$, then it osculates to $Z(f)$ to order four at
$p$ if $f(p)=0$ and
\begin {equation}
\label {flec} F_i(p;v)=0,\quad \text{ for } i=1,2,3,4 .
\end {equation}

The four-dimensional \emph{flecnode polynomial} of $f$, denoted
$\fl{f}^4$, is the polynomial obtained by eliminating $v$ from the
four equations in the system (\ref{flec}). (See
Salmon~\cite{salmon}, and the relevant applications thereof in
\cite{EKS,GK2}, for details concerning flecnode polynomials in three
dimensions; see also Ivey and Landsberg~\cite{IL} for a more modern
generalization of this concept.) Note that these four polynomials
are homogeneous in $v$ (of respective degrees $1$, $2$, $3$, and
$4$). We thus have a system of four equations in eight variables,
which is homogeneous in the four variables $v_0, v_1, v_2, v_3$.
Eliminating those variables results in a single polynomial equation
in $p=(x,y,z,w)$. Using standard techniques, as in Cox et
al.~\cite{CLO2}, the resulting polynomial $\fl{f}^4$ is the
\emph{multipolynomial resultant} $Res_4(F_1, F_2, F_3, F_4)$ of
$F_1,F_2,F_3,F_4,$ regarding these as polynomials in $v$ (where the
coefficients are polynomials in $p$). By definition, $\fl{f}^4$
vanishes at all the flecnodes of $f$. The following results are
immediate consequences of the theory of multipolynomial resultants,
presented in Cox et al.~\cite{CLO2}.

\begin{lemma}\label{le:flec1}
Given a polynomial $f \in \cplx[x,y,z,w]$ of degree $\Deg\ge 4$, its
flecnode polynomial $\fl{f}^4$ has degree $O(\Deg)$.
\end{lemma}


\noindent{\bf Proof.} The polynomial $F_i$, for $i=1,\ldots,4$, is a
homogeneous polynomial in $v$ of degree $d_i=i$ over
$\cplx[x,y,z,w]$.  By \cite[Theorem 4.9]{CLO2}, putting
$d:=\left(\sum_{i=1}^4 d_i \right) - 3 = 7$, the multipolynomial
resultant $\fl{f}^4=Res_4(F_1,F_2,F_3,F_4)$ is equal to $\frac {D_3}
{D_3'}$, where $D_3$ is a polynomial of degree $\binom {d+3}{3} =
\binom {10}{3} = 120$ in the coefficients of the polynomials $F_i$,
and $D_3'$ is a polynomial of degree $d_1 d_2 d_3 + d_1 d_2 d_4 +
d_1 d_3 d_4 + d_2 d_3 d_4= 6 + 8 + 12 + 24=50$ in these coefficients
(see Cox et al.~\cite[Chapter 3.4, exercises 1,3,6,12,19]{CLO2}).
Since each coefficient of any of the polynomials $F_i$ is of degree
at most $\Deg -1$, we deduce that $\fl{f}^4$ is of degree at most
$O(\Deg)$. \proofend

\begin{lemma}\label{le:flec3}
Given a polynomial $f \in \cplx[x,y,z,w]$ of degree $\Deg\ge 4$,
every line that is fully contained in $Z(f)$ is also fully contained
in $Z(\fl{f}^4)$.
\end{lemma}

\noindent{\bf Proof.} Every point on any such line is a flecnode of
$f$, so $\fl{f}^4$ vanishes identically on the line. \proofend


\paragraph{Ruled Surfaces in four dimensions.}
Flecnode polynomials are a major tool for characterizing ruled
surfaces. This is manifested in the following theorem of Landsberg
\cite{Land}, which is a crucial tool for our analysis. It is
established in \cite{Land} as a considerably more general result,
but we formulate here a special instance that suffices for our
needs.

\begin{theorem}[Landsberg \cite{Land}] \label{th:flec2}
Let $f \in \cplx[x,y,z,w]$ be a polynomial of degree $\Deg\ge 4$.
Then $Z(f)$ is ruled by (complex) lines if and only if  $Z(f)
\subseteq Z(\fl{f}^4)$.
\end{theorem}

We note that Theorem~\ref{th:flec2} extends the classical
Cayley--Salmon theorem in three dimensions (see
Salmon~\cite{salmon}). A quick review of this result is given below.
We also note that we will use a refined version of this theorem,
also due to Landsberg, given as Theorem~\ref{anlan} in
Section~\ref{sec:pf1}.

When $f$ is of degree $\le 3$, we have the following simpler
situation.
\begin {lemma}
\label {le:de} For every polynomial $f \in \cplx[x,y,z,w]$ of degree
$\le 3$, $Z(f)$ is ruled by (possibly complex) lines.
\end {lemma}

\noindent{\bf Proof.} Let $v=(v_0,v_1,v_2,v_3)\in \cplx^4$ be a
direction. First notice that for a point $p \in \cplx^4$, the line
through $p$ in direction $v$ is contained in $Z(f)$ if and only if
the first three equations in (\ref{flec}) are satisfied, because all
the other terms in the Taylor expansion of $f(p+tv)$ always vanish
for a polynomial $f$ of degree $\le 3$. This is a system of three
homogeneous polynomials in $v_0,v_1,v_2,v_3$, of degrees $1,2,3$,
respectively. By B\'ezout's theorem, as stated in
Theorem~\ref{th:agbt} below, the number of solutions (complex
projective, counted with multiplicities) of this system is either
six or infinite, so there is at least one (possibly complex) line
that passes through $p$ and is contained in $Z(f)$. \proofend

\paragraph{Back to three dimensions.}
In three dimensions the analysis is somewhat simpler, and goes back
to the 19th century, in Salmon's work~\cite{salmon} ond others. The
flecnode polynomial $\FL_f$ of $f$, defined in an analogous manner,
is of degree $11\deg(f)-24$~\cite{salmon}. Theorem~\ref{th:flec2} is
replaced by the Cayley--Salmon theorem~\cite{salmon}, with the
analogous assertion that $Z(f)$ is ruled by lines if and only if
$Z(f)\subseteq Z(\FL_f)$. A simple proof of the Cayley--Salmon
theorem can be found in Terry Tao's blog~\cite{Tao14}.

We will be using the following result, established by
Guth and Katz~\cite{GK}; see also \cite{EKS}.
\begin{proposition} \label{caysal}
Let $f$ be a trivariate irreducible polynomial of degree $D$. If $Z(f)$ fully
contains more than $11D^2-24D$ lines then $Z(f)$ is ruled by (possibly complex) lines.
\end{proposition}
\noindent{\bf Proof.}
Apply Corollary~\ref{co:btl} to $\FL_f$ and $f$, to conclude that $\FL_f$
and $f$ must have a common factor. Since $f$ is irreducible, this factor must be $f$
itself, and then the Cayley--Salmon theorem implies that $Z(f)$ is ruled.
\proofend

\subsection{Flat points and the second fundamental form} \label{se:algflat}
We continue with the four-dimensional setup. Extending the notation
in Guth and Katz~\cite{GK} (see also \cite{EKS}, and also
Pressley~\cite{Pr} and Ivey and Landsberg~\cite{IL} for more basic
references),  we call a non-singular point $p$ of $Z(f)$
\emph{linearly flat}, if it is incident to at least three distinct
2-flats that are fully contained in $Z(f)$ (and thus also in the
tangent hyperplane $T_p Z(f)$). (The original definition,
in~\cite{EKS, GK2}, for the three-dimensional case, is that a
non-singular point $p\in Z(f)$ is linearly flat if it is incident to
three distinct lines that are fully contained in $Z(f)$) The
condition for a point $p$ to be linearly flat can be worked out as
follows, suitably extending the technique used in three dimensions
in \cite{EKS,GK}. Although this extension is fairly routine, we are
not aware of any previous concrete reference, so we spell out the
details for the sake of completeness.

Let $p$ be a non-singular point of $Z(f)$, and let $f^{(2)}$ denote
the second-order Taylor expansion of $f$ at $p$. That is, we have,
for any direction vector $v$ and $t\in \cplx$,
\begin{equation}
\label{bfie}
\begin {array}{ll}
f^{(2)}(p+tv) &= t\nabla f(p)\cdot v + \tfrac12 t^2 v^T H_{f}(p)v.
\end {array}
\end{equation}
If $p$ is linearly flat, there exist three 2-flats $\pi_1$, $\pi_2$,
$\pi_3$, contained in the tangent hyperplane $T_p Z(f)$, such that
$v^T H_f(p) v=0$, for all $v\in \pi_1,\pi_2,\pi_3$ (clearly, the first term
$\nabla f(p)\cdot v$ also vanishes for any such $v$). Using a suitable
coordinate frame within $T_p Z(f)$, we can regard $v^T H_f(p) v$ as a
quadratic trivariate homogeneous polynomial.

Since $v^T H_f(p) v$ vanishes on three 2-flats inside $T_p Z(f)$, a
(generic) line $\ell$, fully contained in $T_p Z(f)$ and not passing
through $p$, intersects these 2-flats at three distinct points, at
which $v^T H_{f} v$ vanishes. Since this is a quadratic polynomial,
it must vanish identically on $\ell$. Thus, $v^T H_f v$ is zero for
all vectors $v \in T_p Z(f)$, and thus $f^{(2)}$ vanishes
identically on $T_p Z(f)$. In this case, we say that $p$ is a
\emph{flat} point of $Z(f)$. Therefore, every linearly flat point of
$Z(f)$ is also a flat point of $Z(f)$ (albeit not necessarily vice
versa. \footnote{For example, for the surface in $\reals^3$ defined
by the zero set of $f=x+y+z+x^3$, the point $0=(0,0,0)\in Z(f)$ is
flat (because the second order Taylor expansion of $f$ near $0$ is
the plane $x+y+z=0$), but is not linearly flat, since there is no
line incident to $0$ and contained in $Z(f)$.}) The same definition
applies in three dimensions too.

We next express the set of linearly flat points of $Z(f)$ as the
zero set of a certain collection of polynomials. To do so, we define
three canonical 2-flats, on which we test the vanishing of the
quadratic form $v^T H_f v$. (The preceding analysis shows that, for
a linearly flat point, it does not matter which triple of 2-flats is
used for testing the linear flatness, as long as they are distinct.)
These will be the 2-flats
\begin {equation}
\pi_p^x:= T_p Z(f) \cap \{x=x_p\}, \quad \pi_p^y:=T_p Z(f) \cap
\{y=y_p\}, \text{ and } \pi_p^z:=T_p Z(f) \cap \{z=z_p\}.
\end {equation}
These are indeed distinct 2-flats, unless $T_p Z(f)$ is orthogonal
to the $x$-, $y$-, or $z$-axis. Denote by $Z(f)_{axis}$ the subset
of non-singular points $p\in Z(f)$, for which $T_p Z(f)$ is
orthogonal to one of these axes, and assume in what follows that $p
\in Z(f)\setminus Z(f)_{axis}$. We can ignore points in
$Z(f)_{axis}$ by assuming that the coordinate frame of the ambient
space is generic, to ensure that none of our (finitely many) input
points has a tangent hyperplane that is orthogonal to any of the
axes.

\begin {lemma}
\label{le:flat} Let $p$ be a non-singular point of $Z(f)\setminus
Z(f)_{axis}$. Then $p$ is a flat point of $Z(f)$ if and only if $p$
is a flat point of each of the varieties
$Z(f|_{\{x=x_p\}}), Z(f|_{\{y=y_p\}}),Z(f|_{\{z=z_p\}}).$
\end {lemma}

\noindent{\bf Proof.}
Note that the three varieties in the lemma are two-dimensional
varieties within the corresponding three-dimensional cross-sections
$x=x_p$, $y=y_p$, and $z=z_p$, of 4-space.

If $p$ is a flat point of $Z(f)\setminus Z(f)_{axis}$, then the
second-order Taylor expansion $f^{(2)}$ vanishes identically on $T_p
Z(f)$. By the assumption on $p$, we have
\begin {align*}
T_p Z(f|_{\{x=x_p\}}) & = T_p Z(f) \cap \{x=x_p\}, \\
T_p Z(f|_{\{y=y_p\}}) & = T_p Z(f) \cap \{y=y_p\}, \text{ and } \\
T_p Z(f|_{\{z=z_p\}}) & = T_p Z(f) \cap \{z=z_p\},
\end {align*}
and these are three distinct 2-flats. Therefore,
$f|_{\{x=x_p\}}^{(2)}$ vanishes identically on $T_p
Z(f|_{{\{x=x_p\}}})$, implying that $p$ is a flat point of
$Z(f|_{\{x=x_p\}})$; similarly $p$ is a flat point of
$Z(f|_{\{y=y_p\}})$ and of $Z(f|_{\{z=z_p\}})$. For the other
direction, notice that if $p$ satisfies the assumptions in the
lemma, and is a flat point of each of $Z(f|_{\{x=x_p\}}),
Z(f|_{\{y=y_p\}})$, and $Z(f|_{\{z=z_p\}})$, then $f^{(2)}$ vanishes
on three distinct 2-flats contained in $T_p Z(f)$ (namely, the
intersection of $T_p Z(f)$ with $\{x=x_p\}, \{y=y_p\}$ and
$\{z=z_p\}$), which are distinct since $p \not\in Z(f)_{axis}$.
Since $f^{(2)}$ is quadratic, the argument given above implies that
it is identically 0 on $T_p Z(f)$. \proofend


Recall from Elekes et al.~\cite{EKS} that $p$ is flat for
$f|_{\{x=x_p\}}$ if and only if $\Pi_j^1:=\Pi_j(f|_{\{x=x_p\}})$
vanishes at $p$, for $j=1,2,3$, where $\Pi_j(h) = (\nabla h \times
e_j)^T H_h (\nabla h \times e_j)$, and where $e_1,e_2,e_3$ denote
the unit vectors in the respective $y$-, $z$-, and $w$-directions,
and the symbol $\times$ stands for the vector product in
$\{x=x_p\}$, regarded as a copy of $\cplx^3$. In fact, when $x_p$ is
also considered as a variable (call it $x$ then), we get that, as in
the three-dimensional case, each of $\Pi_j^1$, for $j=1,2,3$, is a
polynomial in $x,y,z,w$ of (total) degree $3\Deg-4$. Similarly, the
analogously defined polynomials $\Pi_j^2:=\Pi_j(f|_{\{y=y_p\}}),
\Pi_j^3:=\Pi_j(f|_{\{z=z_p\}})$, for $j=1,2,3$, vanish at $p$ if and
only if $p$ is a flat point of $f|_{\{y=y_p\}}$ and
$f|_{\{z=z_p\}}$. By Lemma \ref{le:flat}, we conclude that a
non-singular point $p\in Z(f)\setminus Z(f)_{axis}$ is flat if and
only if $\Pi_j^i(p)=0,$ for $1 \le i,j \le 3$.

We say that a line $\ell \subset Z(f)$ is a \emph{singular} line of
$Z(f)$, if all of its points are singular. We say that a line $\ell
\subset Z(f)$ is a \emph{flat} line of $Z(f)$ if it is not a singular line of $Z(f)$,
and all of its non-singular points are flat. An easy observation is
that a flat line can contain at most $\Deg-1$ singular points of $Z(f)$
(these are the points on $\ell$ where all four first-order partial derivatives
of $f$ vanish). Similarly, a non-singular line is flat if
(and only if) it is incident to at least $3D-3$ flat points.

\paragraph{The second fundamental form.}
We use the following notations and results from differential
geometry; see Pressley~\cite{Pr} and Ivey and Landsberg~\cite{IL}
for details. For a variety $X$, the differential $d\gamma$ of the
\emph{Gauss mapping} $\gamma$ that maps each point $p \in X$ to its
tangent space $T_p X$, is called the \emph{second fundamental form}
of $X$. In four dimensions, for $X=Z(f)$, and for any non-singular
point $p \in Z(f)$, the second fundamental form, locally near $p$,
can be written as (see~\cite{IL})
$$
\sum_{1\le i,j \le 3} a_{ij} du_idu_j,
$$
where $x=x(u_1,u_2,u_3)$ is a parametrization of $Z(f)$, locally
near $p$, and $a_{ij}=x_{u_iu_j}\cdot \bf n$, where ${\bf n}={\bf
n}(p)=\nabla f(p)/\| \nabla f(p) \|$ is the unit normal to $Z(f)$ at
$p$. Since the second fundamental form is the differential of the
Gauss mapping, it does not depend on the specific local
parametrization of $f$ near $p$. An important property of the
second fundamental form is that it vanishes at every non-singular
flat point $p\in Z(f)$ (see, e.g., Pressley~\cite{Pr} and
Ivey and Landsberg~\cite{IL}).

\begin {lemma}
\label {le:lar} If a line $\ell \subset Z(f)$ is flat, then the
tangent space $T_p Z(f)$ is fixed for all the non-singular points
$p\in\ell$.
\end {lemma}

\noindent{\bf Proof.} The proof applies a fairly standard argument
in differential geometry (see, e.g., Pressley ~\cite{Pr});
see also a proof of a similar claim for the
three-dimensional case in~\cite[Appendix]{EKS}. Fix a non-singular
point $p \in \ell$, and assume that $x=x(u_1,u_2,u_3)$ is a
parametrization of $Z(f)$, locally near $p$. We assume, as we may,
that the relevant neighborhood $N_p$ of $p$ consists only of
non-singular points. For any point $(a,b,c)$ in the corresponding
parameter domain, $x_{u_1},x_{u_2},x_{u_3}$ span the tangent space
to $Z(f)$ at $x(a,b,c)$. Indeed, since $x(u_1,u_2,u_3)$ is a local
parametrization, its differential $(dx)_{(a,b,c)}:
T_{(a,b,c)}\cplx^3 \to T_{x(a,b,c)}Z(f)$ is an isomorphism. Hence,
the image of this latter map is spanned by $x_{u_1},x_{u_2},x_{u_3}$
at $x(a,b,c)$. In particular, we have
$$x_{u_i}\cdot {\bf n} = 0, \quad i=1,2,3,$$
over $N_p$. We now differentiate these equations with respect to
$u_j$, for $j=1,2,3,$ and obtain
$$
x_{u_iu_j}\cdot {\bf n}+x_{u_i}\cdot {\bf n}_{u_j}\equiv 0 \quad
\text{on } \ell\cap N_p,  \text{ for } 1 \le i,j \le 3.
$$
The first term vanishes because $\ell$ is flat, so, as noted above,
the second fundamental form vanishes at each non-singular point of
$\ell$. We therefore have
$$
x_{u_i} \cdot {\bf n}_{u_j}\equiv 0 \quad \text{on } \ell\cap N_p,
\text{ for } i,j=1,2,3.
$$
Since $x_{u_1},x_{u_2},x_{u_3}$ \emph{span} the tangent space $T_q
Z(f)$, for each $q\in N_p$, it follows that ${\bf n}_{u_j}(q)$ is
orthogonal to $T_q Z(f)$ for each $q \in \ell \cap N_p$, and thus
must be parallel to ${\bf n}(q)$ in this neighborhood. However,
since ${\bf n}$ is of unit length, we have ${\bf n}\cdot {\bf
n}\equiv 1$, and differentiating this equation yields
$$
{\bf n}_{u_j}\cdot {\bf n} \equiv 0 \quad \text{on } \ell\cap N_p,
\text{ for } j=1,2,3.
$$
Since ${\bf n}_{u_j}(q)$ is both parallel and orthogonal to ${\bf n}(q)$,
it must be identically zero on $\ell \cap N_p$, for $j=1,2,3$.

Write $\ell = p+t v, \ t\in \cplx$, and define ${\bf h}(t):={\bf n}(p+t v)$,
for $t\in \cplx$. Then, in a suitable tensor notation,
$$
{\bf h}'(t)=({\bf n}_{u_1}(p+tv), {\bf n}_{u_2}(p+tv),
{\bf n}_{u_3}(p+tv))\cdot v \equiv 0,
$$
locally near $t=0$. Thus, ${\bf n}(p+tv)$ is constant locally near $t=0$,
implying that $\bf n$ is constant along $\ell$, locally near $p$.

It still remains to show that $\bf n$ is constant on the set of all
the non-singular points of $Z(f)$ contained in $\ell$. Set
$$
Z_s(\ell):=\{t\in\cplx \mid p+tv \text{ is a singular point of } Z(f)\}.
$$
As $\ell$ is not singular, $|Z_s(\ell)|\le \Deg-1$ (as already
observed). The map $t\mapsto {\bf n}(p+tv)$ is constant in a
neighborhood of every point $t$ of $Z_{ns}(\ell):=\cplx\setminus
Z_s(\ell)$. Since $Z_{ns}(\ell)$ is a connected set,\footnote{
  This property holds for $\cplx$ but not for $\reals$.}
${\bf n}$ has a fixed value at all the non-singular points on
$\ell$, as asserted. Since the tangent hyperplanes $T_p Z(f)$ along
$\ell$ all contain the line $\ell$ itself, and all have the same
normal, we deduce that $T_p Z(f)$ is fixed for all non-singular
points $p\in \ell$. \proofend


\subsection{Finitely and infinitely ruled surfaces in four dimensions, and
u-resultants} \label{se:algres}

Recall again the definition of $\Xi_p$, for a polynomial
$f\in\cplx[x,y,z,w]$ and a point $p\in Z(f)$, which is the union of
all (complex) lines passing through $p$ and fully contained in
$Z(f)$, and that of $\Sigma_p$, as the set of directions (considered
as points in $\mathbb P T_p Z(f)$) of these lines.

Fix a line $\ell\in \Xi_p$, and let $v=(v_0,v_1,v_2,v_3) \in \mathbb P^3$
represent its direction.
Since $\ell\subset Z(f)$, the four terms $F_i(p;v)=\nabla_v^i f(p)$,
for $i=1,2,3,4$, must vanish at $p$.
These terms, which we denote shortly as $F_i(v)$ at the fixed $p$,
are homogeneous polynomials of respective degrees $1,2,3,$ and $4$
in $v=(v_0,v_1,v_2,v_3)$. (Note that when $\Deg\le 3$, some of these
polynomials are identically zero.)

In this subsection we provide a (partial) \emph{algebraic}
characterization of points $p\in Z(f)$ for which $|\Sigma_p|$ is
infinite; that is, points that are incident to infinitely many lines
that are fully contained in $Z(f)$. We refer to this situation by
saying that $Z(f)$ is \emph{infinitely ruled} at $p$. To be precise,
here we only characterize points that are incident to infinitely
many lines that osculate to $Z(f)$ to order three. The passage from
this to the full characterization will be done during the analysis
in the next section.

\paragraph{u-resultants.}
The algebraic tool that we use for this purpose are
\emph{u-resultants}. Specifically, following and specializing
Cox et al.~\cite[Chapter 3.5, page 116]{CLO2}, define, for a vector
$u=(u_0,u_1,u_2,u_3) \in \mathbb P^3$,
$$
U(p;u_0,u_1,u_2,u_3)=Res_4\Bigl(F_1(p;v),F_2(p;v),F_3(p;v),u_0v_0+u_1v_1+u_2v_2+u_3v_3\Bigr)
,
$$
where $Res_4(\cdot)$ denotes, as earlier, the multipolynomial
resultant of the four respective (homogeneous) polynomials, with
respect to the variables $v_0,v_1,v_2,v_3$. For fixed $p$, this is
the so-called \emph{u-resultant} of $F_1(v),F_2(v),F_3(v)$.
\begin {theorem}
\label {ren} The function $U(p;u_0,u_1,u_2,u_3)$ is a homogeneous
polynomial of degree six in the variables $u_0,u_1,u_2,u_3$, and is
a polynomial of degree $O(\Deg)$ in $p=(x,y,z,w)$. For fixed
$p\in Z(f)$, $U(p;u_0,u_1,u_2,u_3)$ is identically
zero as a polynomial in $u_0,u_1,u_2,u_3$, if and only if there are
infinitely many (complex) directions $v=(v_0,v_1,v_2,v_3)$, such
that the corresponding lines $\{p+tv \mid t\in \cplx\}$ osculate to
$Z(f)$ to order three at $p$.
\end {theorem}

\noindent{\bf Proof.} By definition, the osculation property in the
theorem, for given $p$ and $v$, is equivalent to
$F_1(p;v)=F_2(p;v)=F_3(p;v)=0$. Regarding $F_1$, $F_2$, $F_3$ as
homogeneous polynomials in $v$, the degree of $U$ in
$u_0,u_1,u_2,u_3$ is $\deg(F_1)\deg(F_2)\deg(F_3)=3!=6$ (see
Cox et al.~\cite[Exercise 3.4.6.b]{CLO2}). Put
$d=\deg(F_1)+\deg(F_2)+\deg(F_3)+1=7$. Then the total degree of $U$
in the coefficients of $F_i$, each being a polynomial in $p$ of
degree at most $D$, is at most $\binom{d}{3}=\binom 7 3 = 35$ (see
also the proof of Lemma~\ref{le:flec1} and Cox et al.~\cite[Exercises 3.4.6.c,
3.4.19]{CLO2}), and thus the degree of $U$ as a polynomial in $p$ is
$O(\Deg)$.

Put $H(u,v)=u_0v_0+u_1v_1+u_2v_2+u_3v_3$, and, for any $v\in
\cplx^4$, denote by $H_v$ the hyperplane $H(u,v)=0$. Fix $p \in
Z(f)$, and regard $F_1,F_2,F_3, H(u,\cdot)$ as polynomials in $v$.
If the osculation property holds at $p$ (for infinitely many lines)
then $Z(F_1,F_2,F_3)$ is infinite, so it is at least 1-dimensional.
Thus, for any $u=(u_0,u_1,u_2,u_3) \in \cplx^4$, the variety
$Z(F_1,F_2,F_3,H(u,v))$ is non-empty, so the multipolynomial
resultant of these four polynomials (in $v$) vanishes at $u$. Since
this holds for all $u \in \cplx^4$, It follows from Cox et
al.~\cite[Proposition 1.1.5]{CLO2} that $U \equiv 0$.

Suppose then that the osculation property does not hold (for
infinitely many lines) at $p$, so $Z(F_1,F_2,F_3)$ is finite. Pick
any $u_0 \not\in \bigcup_{v \in Z(F_1,F_2,F_3)} H_v$. Then, for
every $v \in Z(F_1,F_2,F_3)$, we have $H(u_0,v) \ne 0$, implying
that
$$
Z(F_1,F_2,F_3,H(u_0,\cdot))=\{v\in Z(F_1,F_2,F_3) \mid H(u_0,v)=0\}=\emptyset.
$$
Therefore, by the properties of multipolynomial resultants,
$U(u_0)\ne 0$, and $U$ is not identically zero. \proofend


\noindent{\bf Remark.} Theorem~\ref{ren} shows that the subset of
$Z(f)$ consisting of the points incident to infinitely many lines
that osculate to $Z(f)$ to order three is contained in a \emph{subvariety}
of $Z(f)$, which is the intersection of $Z(f)$ with the common zero set
of the coefficients of $U$ (considered as polynomials in $x,y,z,w$).

\begin {corollary}
\label {coren} Fix $p\in Z(f)$. The polynomial
$U(p;u_0,u_1,u_2,u_3)$ is identically zero, as a polynomial in
$u_0,u_1,u_2,u_3$, if and only if there are more than six (complex)
lines osculating to $Z(f)$ to order $3$ at $p$.
\end {corollary}

\noindent{\bf Proof.} The polynomial $F_i$ is either $0$ or of
degree $i$ (in $v$, for a fixed value of $p$), for $i=1,2,3$. By
Theorem~\ref{th:agbt}, the number of their common zeros
$v=(v_0,v_1,v_2,v_3)$ is either six (counting complex projective
solutions with multiplicity; see also the proof of
Theorem~\ref{ren}) or infinite. The result then follows from Theorem
\ref{ren}. \proofend


\section{Proof of Theorem~\ref{th:main}} \label{sec:pf1}

Let $P,L,m,n,q$, and $s$ be as in the theorem.

The proof proceeds by induction on $m$, where the base cases of the
induction are the ranges $m\le \sqrt{n}$ and $m\le M_0$, for a
sufficiently large constant $M_0$. In both cases we have $I(P,L) \le
A(m+n)$, for a suitable choice of $A$. \footnote{When $m\le \sqrt n$
(or when $n \le \sqrt m$), an immediate application of the
Szemer\'edi--Trotter theorem yields the linear bound $O(m+n)$.}
Assume then that the bound holds for all $m' < m$, and consider an
instance involving sets $P$, $L$, with $|P|=m
> \sqrt{|L|} = \sqrt{n}$, and $m> M_0$.

As already discussed, the bound in (\ref{ma:in}) is qualitatively
different in the two ranges $m = O(n^{4/3})$ and $m =
\Omega(n^{4/3})$, and the analysis will occasionally have to
bifurcate accordingly. Nevertheless, the bifurcation is mainly in
the choice of various parameters, and in manipulating them.
Most of the technical
details that deal with the algebraic structure of the problem are
identical. We will therefore present the analysis jointly for both
cases, and bifurcate only locally, when the induction itself, or
tools that prepare for the induction, get into action, and require
different treatments in the two cases.

As promised in the overview, we will use two different partitioning
schemes, one with a polynomial of ``large'' degree, and one with a
polynomial of ``small'' degree. We start naturally with the
first scheme.

An important issue to bear in mind is that, unlike most of the
material in the preceding section, where the underlying field was
$\cplx$, the analysis in this section is over the reals.
Nevertheless, this is essentially needed only for constructing a polynomial
partitioning, which is meaningless over $\cplx$. Once this is done, the
analysis of incidences between points and lines on the zero set of
the partitioning polynomial can be carried out over the complex field
just as well as over $\reals$, and then the machinery reviewed and
developed in the previous section can be brought to bear.

\paragraph{First partitioning scheme.}
Fix a parameter $r$, given by
$$
r = \begin{cases}
cm^{8/5}/n^{4/5} & \text{if $m \le an^{4/3}$} \\
cn^4/m^2  & \text{if $m \ge an^{4/3}$} ,
\end{cases}
$$
where $a$ and $c$ are suitable constants. Note that, in both cases,
$1\le r\le m$, for a suitable choice of the constants of
proportionality, unless either $m = \Omega(n^2)$ or $n =
\Omega(m^2)$, extreme cases that have already been handled. We refer
to the cases $n^{1/2} \le m \le an^{4/3}$ and $an^{4/3} < m \le n^2$
as the cases of \emph{small $m$} and of \emph{large $m$},
respectively.

We now apply the polynomial partitioning theorem of Guth and Katz (see \cite{GK2}
and \cite[Theorem 2.6]{KMS}), to obtain an $r$-partitioning 4-variate (real)
polynomial $f$ of degree
\begin {equation}
\label {eq:degofpart}
\begin {array}{ll}
\Deg=O(r^{1/4}) \le \begin{cases}
c_0m^{2/5}/n^{1/5} & \text{if $m \le an^{4/3}$} \\
c_0n/m^{1/2}  & \text{if $m \ge an^{4/3}$} ,
\end{cases}
\end {array}
\end {equation}
for another suitable constant $c_0$. That is, every connected
component of $\reals^4\setminus Z(f)$ contains at most $m/r$ points
of $P$, where, as above, $Z(f)$ denotes the zero set of $f$. By
Warren's theorem~\cite{w-lbanm-68} (see also \cite{KMS}), the number
of components of $\reals^4\setminus Z(f)$ is $O(\Deg^4) = O(r)$.

Set $P_0:= P\cap Z(f)$ and $P':=P\setminus P_0$. We recall that,
although the points of $P'$ are more or less evenly partitioned
among the cells of the partition, no nontrivial bound can be
provided for the size of $P_0$; in the worst case, all the points of
$P$ could lie in $Z(f)$.  Each line $\ell\in L$ is either fully
contained in $Z(f)$ or intersects it in at most $\Deg$ points (since
the restriction of $f$ to $\ell$ is a univariate polynomial of
degree at most $\Deg$). Let $L_0$ denote the subset of lines of $L$
that are fully contained in $Z(f)$ and put $L' = L\setminus L_0$.

We have
\begin{equation} \label{eqsplit}
I(P,L) = I(P_0,L_0) + I(P_0,L') + I(P',L') .
\end{equation}
As can be expected (and noted earlier), the harder part of the
analysis is the estimation of $I(P_0,L_0)$. Indeed, it might happen
that $Z(f)$ is a hyperplane, and then the best (and worst-case
tight) bound we can offer is the bound specified by Theorem
\ref{ttt}. It might also happen that $Z(f)$ contains some 2-flat, in
which case we are back in the planar scenario, for which the best
(and worst-case tight) bound we can offer is the
Szemer\'edi--Trotter bound (\ref{inc2}). Of course, the assumptions
of the theorem come to the rescue, and we will see below how exactly
they are used.

We first bound the second and third terms of (\ref{eqsplit}). We
have
\begin {equation}
\label {eq:meq} I(P_0,L') \le |L'|\cdot \Deg \le n\Deg,
\end {equation}
because, as just noted, a line not fully contained in $Z(f)$ can
intersect this set in at most $\Deg$ points. To estimate $I(P',L')$,
we put, for each cell $\tau$ of the partition, $P_\tau = P\cap \tau$,
and let $L_\tau$ denote the set of the lines of
$L'$ that cross $\tau$; put $m_\tau = |P_\tau| \le m/r$, and $n_\tau
= |L_\tau|$. Since every line $\ell\in L'$ crosses at most $1+\Deg$
components of $\reals^4\setminus Z(f)$ (because it has to pass
through $Z(f)$ in between cells), we have
\begin {equation}
\label {aeq} \sum_\tau n_\tau \le |L'|(1+\Deg) \le n(1+\Deg) .
\end {equation}
Clearly, we have
\begin {equation*}
I(P',L') = \sum_\tau I(P_\tau,L_\tau).
\end {equation*}
We now bifurcate depending on the value of $m$.

\paragraph{Estimating $I(P',L')$: The case of small $m$.}
Here we use the easy upper bound (which holds for any pair of sets
$P_\tau,L_\tau$)
$$
I(P_\tau,L_\tau) = O(|P_\tau|^2+|L_\tau|) = O((m/r)^2 + n_\tau) .
$$
Summing these bounds over the cells, using (\ref{aeq}), and
recalling the value of $r$ (and of $D$), we get
$$
I(P',L') = \sum_\tau I(P_\tau,L_\tau) = O(m^2/r + nr^{1/4}) =
O(m^{2/5}n^{4/5}) .
$$
\paragraph{Estimating $I(P',L')$: The case of large $m$.}
Here we use the dual (generally applicable) upper bound
$I(P_\tau,L_\tau) = O(|L_\tau|^2+|P_\tau|)$, which, by splitting
$L_\tau$ into subsets of size at most $|P_\tau|^{1/2}$, becomes
$$
I(P_\tau,L_\tau) = O(|P_\tau|^{1/2}|L_\tau| + |P_\tau|) =
O((m/r)^{1/2}n_\tau + m_\tau) .
$$
Summing these bounds over the cells, using (\ref{aeq}), and recalling the
value of $r$, we get
$$
I(P',L') = \sum_\tau I(P_\tau,L_\tau) = O((m/r)^{1/2}nr^{1/4} + m) =
O(m^{1/2}n/D + m) = O(m) .
$$
Combining both bounds, we have:
\begin{equation} \label{incells}
I(P_0,L') + I(P',L') = O\left( m^{2/5}n^{4/5} + m \right) .
\end{equation}
Note that in this part of the analysis we do not need the
assumptions involving $q$ and $s$ --- the large degree trivializes
the analysis within the cells of the partition.


\paragraph{Estimating $I(P_0,L_0)$.}
We next bound the number of incidences between points and lines that
are contained in $Z(f)$. To simplify the notation, write $P$ for $P_0$
and $L$ for $L_0$, and denote their respective cardinalities as $m$ and $n$.
(The reader should keep this convention in mind, as we will ``undo'' it towards the end of the analysis.)
To be precise, we will not be able to account explicitly for all types of these incidences
(for the present choices of $D$). Our strategy is to
obtain an explicit bound for a subset of the incidences, which is
subsumed by the bound in (\ref{ma:in}), and
then prune away those lines and points that participate in these
incidences. We will be left with ``problematic'' subsets of points
and lines, and we will then handle them in a second, new, induction-based
partitioning step. A major goal for the first stage is to show that,
for the set of surviving lines, the parameters $q$ and $s$
can be replaced by the respective parameters $O(D^2)$ and
$O(D)$ that ``pass well'' through the induction; see below for details.\footnote{
  Note that in general the bounds $O(D^2)$ and $O(D)$ are not necessarily
  smaller than their respective original counterparts $q$ and $s$. Nevertheless,
  they uniformly depend on $m$ and $n$ in a way that makes them fit the induction
  process, whereas the parameters $q$ and $s$, over which we have no control, do not.}

By the nature of its construction, $f$ is in general reducible (see
\cite{GK2}). However, to apply successfully certain steps of the
forthcoming analysis, we will need to assume that $f$ is
irreducible, so we will apply the analysis separately to each
irreducible factor of $f$, and then sum up the resulting bounds.
(The actual problem decomposition is subtler --- see below.)

Write the irreducible factors of $f$, in an arbitrary order, as
$f_1,\ldots,f_k$, for some $k \le \Deg$. The points of $P$ are
partitioned among the zero sets of these factors, by assigning each
point $p\in P$ to the first factor in this order whose zero set
contains $p$. A line $\ell\in L$ is similarly assigned to the first
factor whose zero set fully contains $\ell$ (there always exists
such a factor). Then $I(P,L)$ is the sum, over $i=1,\ldots,k$, of
the number of incidences between the points and the lines that are
assigned to the (same) $i$th factor, plus the number of incidences
between points and lines assigned to different factors. The latter
kind of incidences is easier to handle. Indeed, if $(p,\ell)$ is an
incident pair in $P \times L$, so that $p$ is assigned to $f_i$ and
$\ell$ is assigned to $f_j$, for $i \ne j$ (necessarily $i<j$), then
the incidence occurs at an intersection of $\ell$ with $Z(f_i)$.
By construction, $\ell$ is not fully contained in $Z(f_i)$, so it
intersects it in at most $\deg (f_i)$ points, so the overall number
of incidences on $\ell$ of this kind is at most
$\sum_{i \ne j} \deg f_i < \Deg$, and the overall number of
such incidences is therefore at most $n \Deg$.

For the former kind of incidences, we assume in what follows that we
have a single irreducible polynomial $f$, and denote by $P$ and $m$,
for short, the set of points assigned to $f$ and its cardinality,
and by $L$ and $n$ the set of lines assigned to $f$ (and thus fully
contained in $Z(f)$) and its cardinality. We continue to denote the
degree of $f$ as $\Deg$.
(Again, we will undo these conventions towards the end of the analysis.)

This is not yet the end of the reduction, because, in most of the
analysis about to unfold, we need to assume that the points of $P$
are non-singular points of $Z(f)$. To reduce the setup to this
situation we proceed as follows. We construct a sequence of partial
derivatives of $f$ that are not identically zero on $Z(f)$. For this
we assume, as we may, that $f$, and each of its derivatives, are
square-free; whenever this fail, we replace the corresponding
derivative by its square-free counterpart before continuing to
differentiate.
Without loss of generality, assume that this sequence is
$f,f_x,f_{xx}$, and so on.  Denote the $j$-th element in this
sequence as $f_j$, for $j=0,1,\ldots$ (so $f_0=f$, $f_1=f_x$, and so
on). Assign each point $p\in P$ to the first polynomial $f_j$ in the
sequence for which $p$ is non-singular; more precisely, we assign
$p$ to the first $f_j$ for which $f_j(p)=0$ but $f_{j+1}(p)\ne 0$
(recall that $f_0(p)$ is always $0$ by assumption. Similarly, assign
each line $\ell$ to the first polynomial $f_j$ in the sequence for
which $\ell$ is fully contained in $Z(f_j)$ but not fully contained
in $Z(f_{j+1})$ (again, by assumption, there always exists such
$a_j$). If $\ell$ is assigned to $f_j$ then it can only contain
points $p$ that were assigned to some $f_k$ with $k\ge j$. Indeed,
if $\ell$ contained a point $p$ assigned to $f_k$ with $k<j$ then
$f_{k+1}(p)\ne 0$ but $\ell$ is fully contained in $Z(f_{k+1})$,
since $k+1\le j$; this is a contradiction that establishes the
claim.

Fix a line $\ell\in L$, which is assigned to some $f_j$. An incidence
between $\ell$ and a point $p\in P$, assigned to some $f_k$, for $k>j$,
can be charged to the intersection of $\ell$ with $Z(f_{j+1})$ at $p$
(by construction, $p$ belongs to $Z(f_{j+1})$).
The number of such intersections is at most
$D-j-1$, so the overall number of incidences of this sort, over all lines
$\ell\in L$, is $O(nD)$. It therefore suffices to consider only incidences
between points and lines that are assigned to the same zero set $Z(f_i)$.

The reductions so far have produced a finite collection of up to
$O(D)$ polynomials, each of degree at most $D$, so that the points
of $P$ are \emph{partitioned} among the polynomials and so are the
lines of $L$, and we only need to bound the number of incidences
between points and lines assigned to the \emph{same} polynomial.
This is not the end yet, because the various partial derivatives
might be reducible, which we want to avoid. Thus, in a final
decomposition step, we split each derivative polynomial $f_j$ into
its irreducible factors, and reassign the points and lines that were
assigned to $Z(f_j)$ to the various factors, by the same ``first
come first served'' rule used above. The overall number of
incidences that are lost in this process is again $O(nD)$. The
overall number of polynomials is $O(D^2)$, as can easily be checked.
Note also that the last decomposition step preserves non-singularity
of the points in the special sense defined above; that is, as is
easily verified, a point $p\in Z(f_j)$ with $f_{j+1}(p)\ne 0$,
continues to be a non-singular point of the irreducible component it
is reassigned to.

We now fix one such final polynomial, still call it $f$, denote its degree
by $D$ (which is upper bounded by the original degree $D$),
and denote by $P$ and $L$ the subsets of the original sets of
points and lines that are assigned to $f$, and by $m$ and $n$ their
respective cardinalities. (Again, this simplifying convention will be undone
towards the end of the analysis.) We now may assume that $P$ consists exclusively
of \emph{non-singular} points of the \emph{irreducible} variety $Z(f)$.


If $\Deg \le 3$, then, by Lemma \ref{le:de}, $Z(f)$ is ruled by
lines. Hypersurfaces ruled by lines will be handled in the
later part of the analysis. (Note that the cases $\Deg=1$ or
$\Deg=2$ can be controlled by assumption (i') of the theorem (see
below), whereas the case $\Deg=3$ requires a different treatment.)
Suppose then that $\Deg \ge 4$. The flecnode polynomial $\fl{f}^4$
of $f$ (see Section~\ref{se:algruled}) vanishes identically on every
line of $L$ (and thus also on $P$, assuming that each point of $P$
is incident to at least one line of $L$). If $\fl{f}^4$ does not
vanish identically on $Z(f)$, then $Z(f,\fl{f}^4) := Z(f)\cap
Z(\fl{f}^4)$ is a two-dimensional variety (see, e.g.,
Hartshorne~\cite[Exercise I.1.8]{Hart83}). It contains $P$ and all
the lines of $L$ (by Lemma \ref{le:flec3}), and is of degree
$O(\Deg^2)$ (by Theorem~\ref{th:agbt}).
The other possibility is that $\fl{f}^4$
vanishes identically on $Z(f)$, and then Theorem \ref{th:flec2}
implies that $Z(f)$ is ruled by lines. This
latter case, which requires several more refined tools from
algebraic geometry, will be analyzed later.

\subsection*{First case: $Z(f,\fl{f}^4)$ is two-dimensional}

Put $g=\fl{f}^4$. In the analysis below, we only use the facts that
$\deg(g)=O(\Deg)$, and that $Z(f,g)$ is two-dimensional, so the
analysis applies for any such $g$; this comment will be useful in
later steps of the analysis. Recall that in this part of the analysis
$f$ is assumed to be an irreducible polynomial of degree $\ge 4$.

We have a set $P$ of $m$ points and a set $L$ of $n$ lines in $\cplx^4$,
so that $P$ is contained in the two-dimensional algebraic
variety $Z(f,g)\subset \cplx^4$. By pruning away all the lines
containing at most $\max(\Deg,\deg(g))$ points of $P$, we lose
$O(n\Deg)$ incidences, and all the surviving lines are contained in
$Z(f,g)$, as is easily checked. For simplicity of notation, we
continue to denote by $L$ the set of surviving lines.

Let $Z(f,g) = \bigcup_{i=1}^s V_i$ be the decomposition of $Z(f,g)$
into its irreducible components, as described in
Section~\ref{se:algbez}. By Theorem \ref{th:agbt}, we have
$\sum_{i=1}^s \deg(V_i) \le \deg(f)\deg(g)=O(D^2)$.
%

\paragraph{Incidences within non-planar components of $Z(f,g)$.}
Our next step is to analyze the number of incidences between points
and lines within the components of $Z(f,g)$ that are not 2-flats.
For this we first need the following bound on point-line incidences
within a two-dimensional surface in three dimensions. This part of
the analysis is taken from our companion paper~\cite{SS3dv}. We also
refer to Section~\ref{se:algruled} for properties of ruled surfaces.

For a point $p$ on an irreducible singly ruled surface $V$, which is
not the exceptional point of $V$, we let $\Lambda_V(p)$ denote the
number of generator lines passing through $p$ and fully contained in
$V$ (so if $p$ is incident to an exceptional line, we do not count
that line in $\Lambda_V(p)$). We also put $\Lambda_V^*(p) :=
\max\{0,\Lambda_V(p)-1\}$. Finally, if $V$ is a cone and $p_V$ is
its exceptional point (that is, apex), we put $\Lambda_V(p_V) =
\Lambda_V^*(p_V):=0$. We also consider a variant of this notation,
where we are also given a finite set $L$ of lines (where not all
lines of $L$ are necessarily contained in $V$), which does not
contain any of the (at most two) exceptional lines of $V$.  For a
point $p\in V$, we let $\lambda_V(p;L)$ denote the number of lines
in $L$ that pass through $p$ and are fully contained in $V$, with
the same provisions as above, namely that we do not count incidences
with exceptional lines, nor do we cound incidences with an
exceptional point, and put $\lambda_V^*(p;L) :=
\max\{0,\lambda_V(p;L)-1\}$. If $V$ is a cone with apex $p_V$, we
put $\lambda_V(p_V;L) = \lambda^*_V(p_V;L) = 0$. We clearly have
$\lambda_V(p;L) \le \Lambda_V(p)$ and $\lambda^*_V(p;L) \le
\Lambda^*_V(p)$, for each point $p$.

\begin{lemma}
\label{firstflip} Let $V$ be an irreducible singly ruled two-dimensional
surface of degree $\Deg>1$ in $\reals^3$ or in $\cplx^3$. Then, for
any line $\ell$, except for the (at most) two exceptional lines of
$V$, we have
\begin{align*}
& \sum_{p \in \ell \cap V} \Lambda_V(p) \le \Deg \quad\quad\text{if $\ell$ is not fully contained in $V$} , \\
& \sum_{p \in \ell \cap V} \Lambda^*_V(p) \le \Deg
\quad\quad\text{if $\ell$ is fully contained in $V$} .
\end{align*}
\end{lemma}

The following lemma provides the needed infrastructure for our
analysis, and is taken from Sharir and Solomon~\cite[Theorem
15]{SS3dv}.
\begin{lemma}
\label{salta} Let $V$ be a possibly reducible two-dimensional
algebraic surface of degree $D>1$ in $\reals^3$ or in $\cplx^3$,
with no linear components. Let $P$ be a set of $m$ distinct points on $V$ and
let $L$ be a set of $n$ distinct lines fully contained in $V$. Then there
exists a subset $L_0\subseteq L$ of at most $O(\Deg^2)$ lines, such
that the number of incidences between $P$ and $L\setminus L_0$ satisfies
\begin{equation} \label{lmlstar}
I(P,L\setminus L_0) = O\left(m^{1/2}n^{1/2}D^{1/2} + m + n\right) .
\end{equation}
\end{lemma}

\noindent{\bf Sketch of Proof.}
We provide the following sketch of the proof; the full details are given
in the companion paper \cite{SS3dv}. Consider the irreducible components
$W_1,\ldots,W_u$ of $V$. We first argue that the number of lines
that are either contained in the union of the non-ruled components, or
those contained in more than one ruled component of $V$ is $O(\Deg^2)$,
and we place all these lines, as well as the exceptional lines of any
singly ruled component, in the exceptional set $L_0$. We may thus assume
that each surviving line in $L_1:=L\setminus L_0$ is contained in a
unique ruled component of $V$, and is a generator of that component.

The strategy of the proof is to consider each line $\ell$ of $L_1$,
and to estimate the number of its incidences with the points of $P$ in
an indirect manner, via Lemma~\ref{firstflip}, applied to $\ell$ and
to each of the ruled components $W_j$ of $V$.

Specifically, we fix some threshold parameter $\xi$, and dispose of
points that are incident to at most $\xi$ lines of $L_1$, losing at most
$m\xi$ incidences. Let $P_1$ denote the set of surviving points.

Now if a line $\ell\in L_1$ is incident to a point $p\in P_1$, it meets
at least $\xi$ other lines of $L_1$ at $p$. It follows from Lemma~\ref{firstflip}
that the overall number of such lines, over all points in $P_1\cap\ell$, is
roughly $D$, so the number of such points on $\ell$ is at most roughly $D/\xi$,
for a total of $nD/\xi$ incidences of this kind. Choosing $\xi=(nD/m)^{1/2}$
yields the bound $O(m^{1/2}n^{1/2}D^{1/2})$, and the lemma follows.
\proofend

We can now proceed, by deriving two upper bounds for certain types of incidences
between $P$ and $L$. The first bound is relevant for the range $m= O(n^{4/3})$,
and the second bound is relevant for the range $m= \Omega(n^{4/3})$. Nevertheless,
both bounds apply to the entire range of $m$ and $n$.

\begin {proposition}
\label{th:degsquared} The number of incidences involving
non-singular points of $Z(f)$ that are contained in components of
$Z(f,g)$ that are not 2-flats is
\begin {equation}
\min\left\{ O(m\Deg^2+n\Deg),\; O(m+n\Deg^4) \right\} .
\end {equation}
\end {proposition}

\noindent{\bf Proof.} We first establish the bound $O(mD^2+nD)$. Let
$p\in Z(f)$ be a non-singular point. The irreducible decomposition
of $S_p:=Z(f,g)\cap T_p Z(f)$ is the union of one- and
two-dimensional components. Clearly, $S_p$ contains all the lines
that are incident to $p$ and are fully contained in $Z(f,g)$; it is
a variety, embedded in $3$-space (namely, in $T_p Z(f)$), of degree
$O(\Deg^2)$. The union of the one-dimensional components is a curve
of degree $O(\Deg^2)$, so, by Lemma~\ref{le:btl}, it can contain at
most $O(\Deg^2)$ lines; when summing over all $p\in P$, the total
number of incidences with those lines is $O(m\Deg^2)$.

It remains to bound incidences involving the two-dimensional
components of $S_p$ that are not 2-flats. By Sharir and
Solomon~\cite[Lemma 5]{SS3d}, the number of lines incident to $p$
inside these two-dimensional components of $S_p$ is at most
$O(\Deg^2)$, except possibly for lines that lie in a component that
is a cone and has $p$ as its apex. Summing over all $p\in P$, we get
a total of $O(m\Deg^2)$ incidences for this case too, ignoring lines
that lie only in conic (or flat) components.

Note that each two-dimensional component of $S_p$ is necessarily
also a two-dimensional irreducible component of $Z(f,g)$. Hence the
analysis performed so far takes care of all incidences except for
those that occur on conic two-dimensional components of $Z(f,g)$
(and on 2-flats, which we totally ignore in this proposition). Let
$V$ be a conic component of $Z(f,g)$ with apex $p_V$, which is not a
2-flat. We note that $V$ cannot fully contain a line that is not
incident to $p_V$. Indeed, suppose to the contrary that $V$
contained such a line $\ell$. Since $V$ is a cone with apex $p_V$,
for each point $a \in \ell$, the line connecting $a$ to $p_V$ is
fully contained in $V$, and therefore the 2-flat containing $p_V$
and $\ell$ is fully contained in $V$. As $V$ is irreducible and is
not a 2-flat, we obtain a contradiction, showing that no such line
exists. We conclude that any point on $V$, except for $p_V$, is
incident to at most one line that is fully contained in $V$ (a
``generator'' through $p_V$), for a total of $O(m)$ incidences.
Since $Z(f,g)$ is of degree $O(\Deg^2)$, the number of conic
components of $Z(f,g)$ is $O(\Deg^2)$, so, summing this bound over
all components $V$, we get again the bound $O(m\Deg^2)$ on the
number of relevant ``non-apex'' incidences.

Therefore, it remains to bound the number of incidences between the points of
$$
P_c:= \{p_V \mid p_V \text{ is an apex of an irreducible conic
component } V \text{ of } Z(f,g)\}
$$
and the lines of $L$. Since there are at most
$O(\Deg^2)$ irreducible components of $Z(f,g)$, we have $|P_c| \le
c\Deg^2$, for some suitable constant $c$. We next let $L_c$ denote
the set of lines in $L$ containing fewer than $c\Deg$ points of
$P_c$, and claim that any point $p\in P_c$ is incident to fewer than $\Deg$
lines of $L\setminus L_c$. Indeed, otherwise, we would get at least
$\Deg$ lines incident to $p$, each containing at least $c\Deg+1$
points of $P_c$, i.e., at least $c\Deg$ points other than $p$. As
these points are all distinct, we would get that $|P_c| \ge
1+\Deg\cdot c\Deg > c\Deg^2$, a contradiction. On the other hand, by
definition of $L_c$, we have
$$
I(P_c,L_c) = O(n\Deg).
$$
We have thus shown that the number of incidences involving points of
$P_c$ is
$$I(P_c,L)=I(P_c,L_c)+I(P_c,L\setminus L_c)=O(n\Deg)+O(m\Deg)=O(n\Deg+m\Deg),$$
well within the bound that we seek to establish.

\noindent{\bf The second bound.} We next establish the second bound
$O(m+nD^4)$. Let $V$ be an irreducible two-dimensional component of
$Z(f,g)$. If $V$ is not ruled, then by Proposition~\ref{caysal}, it
contains at most $11\deg(V)^2-24\deg(V) < 11\deg(V)^2$ lines.
Summing over all irreducible components of $Z(f,g)$ that are not
ruled, we get at most $11\sum_V \deg(V)^2 = O(D^4)$ lines. Let
$\ell$ be one of those lines, and let $p\in \ell\cap P$. For any
other line $\lambda\in L$ that passes through $p$, we charge its
incidence with $p$ to its intersection with $\ell$. This yields a
total of $O(nD^4)$ incidences, to which we add $O(m)$ for incidences
with those points that lie on only one line of $L$, for a total of
$O(m+nD^4)$ incidences.

We next analyze the irreducible components of $Z(f,g)$ that are
ruled but are not 2-flats. Let $V_1,\ldots, V_k$ denote these
components, for some $k = O(D^2)$. Project all these components onto
some generic hyperplane, and regard them as a single (reducible)
ruled surface in 3-space, whose degree is $\sum_{i=1}^k \deg(V_i) =
O(D^2)$. Lemma~\ref{salta} then yields a subset $L_0$ of $L$ of size
$O(D^4)$, and shows that
$$
I(P,L\setminus L_0) = O\left( m^{1/2}n^{1/2}D + m + n \right).
$$
The lines of $L_0$ are simply added to the set of $O(D^4)$ lines not belonging
to ruled components. This does not affect the asymptotic bound $O(nD^4)$ derived
above. In total we get
$$
O\left( m^{1/2}n^{1/2}D + m + nD^4 \right)
$$
incidences. Since
$$
m^{1/2}n^{1/2}D \le \frac12 \left( m + nD^2 \right),
$$
we obtain the second bound asserted in the proposition. \proofend

\noindent{\bf Remark.} The term $O(nD^4)$ appears to be too weak,
and can probably be improved, using ideas similar to those in the
proof of Lemma~\ref{salta}. Since such an improvement does not have
a significant effect on our analysis, we leave it as an interesting
problem for further research.

\paragraph{Restrictedness of hyperplanes and quadrics, and lines on 2-flats.}
The bounds in Proposition~\ref{th:degsquared}
might be too large, for the current choices of $D$, because of the respective terms
$O(mD^2)$ and $O(nD^4)$. (Technically, the $m$ and $n$ in the definition of $D$
are not necessarily the same as the $m$ and $n$ that denote the size of the current
subsets of the original $P$ and $L$, but let us assume that they are the
same for the present discussion.)
For example, when $m=O(n^{4/3})$ and $D=\Theta(m^{2/5}/n^{1/5})$ (recall that
this is the ``large'' value of $D$ for this range),
we have $mD^2 = \Theta(m^{9/5}/n^{2/5})$, and this is $\gg m^{2/5}n^{4/5}$
when $m\gg n^{6/7}$. Similarly,
when $m=\Omega(n^{4/3})$ and $D=\Theta(n/m^{1/2})$ (which is the value chosen for this range),
we have $nD^4 = \Theta(n^{5}/m^{2})$, and this is $\gg m$ when $m\ll n^{5/3}$.
These bounds will be used in the second partitioning
step, where we use a smaller-degree partitioning polynomial, and for $m$ outside
the problematic ranges, i.e., for $m\le n^{6/7}$ or $m\ge n^{5/3}$; see below for details.
Otherwise, for the current $D$, these bounds need to be finessed and replaced by the following
alternative analysis.\footnote{
  As the calculations worked out above indicate, the bounds
  in Proposition~\ref{th:degsquared} will be within
  the bound (\ref{ma:in}) when $m$ is sufficiently small (below $n^{4/3}$) or sufficiently large (above $n^{5/3}$). For such values of $m$
  we can bypass the induction process, and obtain the desired bounds directly,
  in a single step. See a more detailed description towards the end of this section.}

In the first step of this analysis, we estimate the number of lines
contained in a hyperplane or a quadric (when $Z(f,g)$ is two-dimensional),
and establish the following properties.



\begin {lemma}
\label {le:exhq} Each hyperplane or quadric $H$ is
$O(\Deg^2)$-restricted for the lines of $L$ that are contained in
non-planar components of $Z(f,g)$.
\end {lemma}

\noindent{\bf Proof.} Fix a hyperplane or quadric $H$. Recall that
all the lines in the current set $L$ are contained in $Z(f,g)$. Let
$V$ be an irreducible component of $Z(f,g)$, which is not a 2-flat.
If $V\cap H$ is a curve, then (recalling Theorem~\ref{th:agbt}) its
degree is at most $\deg(V)$ (when $H$ is a hyperplane) or $2\deg(V)$
(when $H$ is a quadric), and can therefore contain at most
$2\deg(V)$ lines, by Lemma~\ref{le:btl}. Therefore, the union of all
the irreducible components $V$ of $Z(f,g)$ which intersect $H$ in a
curve, contains at most $2\sum_V \deg(V) = O(\Deg^2)$ lines. Assume
then that $V\cap H$ is two-dimensional. Since $V$ is irreducible, we
must have $V\cap H=V$, so $V$ is fully contained in $H$. Moreover,
$V$ is an irreducible two-dimensional surface contained in $Z(f)\cap
H$, and therefore must be an irreducible component of $Z(f)\cap H$,
which is a two-dimensional surface of degree $\le \Deg$. By
Theorem~\ref{th:agbt}, $\sum_{V\subset H} \deg(V)\le \deg(Z(f)\cap
H)\le\deg(f)\le\Deg$. If $V$ is not ruled by lines (and, by
assumption, is not a 2-flat), then by Proposition~\ref{caysal}, it
contains at most $11\deg(V)^2$ lines, and summing over all such
components $V$ within $H$, we get a total of at most $\sum_V
11\deg(V)^2 = O(\Deg^2)$ lines.

The remaining irreducible (two-dimensional) components $V$ of
$Z(f,g)$ that meet $H$ (if such components exist) are fully
contained in $H$, and are ruled by lines. As already observed, these
components are also irreducible components of $Z(f)\cap H$, and so,
with the exception of $O(\Deg^2)$ lines (those contained in the
components already analyzed), all the lines of $L$ that lie in $H$
are contained in components of $Z(f)\cap H$ that are ruled by lines.
Since $f$ restricted to $H$ is a polynomial of degree $\le D$, and
since we are interested in lines of $L$ that are not contained in
planar components of $Z(f)\cap H$ we conclude that $H$ is
$O(\Deg^2)$-restricted, with respect to the subset of $L$ mentioned
in the lemma. \proofend


We next analyze the number of lines contained in a 2-flat.

\begin {lemma}
\label {le:line2flat}
Let $\pi$ be a 2-flat that is not fully contained in $Z(f,g)$.
Then the number of lines fully contained in $Z(f)\cap\pi$ is $O(\Deg)$.
\end {lemma}

\noindent{\bf Proof.} The intersection $Z(f)\cap \pi$ is either
$\pi$ itself, or a curve of degree $\le\Deg$. The latter case
implies (using Lemma~\ref{le:btl}) that $\pi$ contains at most
$\Deg$ lines that are fully contained in $Z(f)$. In the former case
$\pi \subset Z(f)$. By assumption, $\pi$ is not contained in
$Z(f,g)$, implying that $g$ intersects $\pi$ in a curve of degree
$O(\Deg)$ (since $\pi \cap Z(f,g)=\pi \cap Z(g)$), and can therefore
contain at most $O(\Deg)$ lines that are fully contained in $Z(f)$.
\proofend

\paragraph{Recap.}
Summing up what has been done so far, we can classify the incidences
in $I(P,L)$ into the following types. Recall that the analysis is
confined to a single irreducible factor $f$ of the original
polynomial or of some partial derivative of such a factor.

\noindent
(a) We treat the cases where $f$ is linear or quadratic separately,
using a variant of Theorem~\ref{ttt}, which takes into account the
restrictedness of hyperplanes and quadrics; see Proposition~\ref{th:bwhyp}
below.

\noindent
(b) We treat the case where $Z(f)$ is ruled by lines separately
(this is the second case in the analysis, when $Z(f,\flf)$ is three-dimensional).

If $f$ is not ruled by lines and is of degree $\ge 4$
(recall that each surface of degree at most $3$ is ruled by lines---see Lemma~\ref{le:de}),
then there are two kinds of incidences that need to be considered.

\noindent (c) Incidences between points and lines that are contained
in irreducible components of $Z(f,\flf)$ (or, more generally, of
$Z(f,g)$, for other suitable polynomials $g$) that are not 2-flats.
We have bounded the number of these incidences in
Proposition~\ref{th:degsquared} in two different ways, but we also
ignored these incidences, passing them to the induction in the
second partitioning step, to be presented later, where we now know
that each hyperplane and quadric is $O(D^2)$-restricted, and each
2-flat contains at most $O(D)$ lines of $L$. For both properties to
hold, we first have to get rid of all the lines of $L$ that are
contained in 2-flats within $Z(f,g)$, and we will perform this
pruning after bounding the number of incidences involving lines that
are contained in such 2-flats. This will make the
$O(D^2)$-restrictedness in Lemma~\ref{le:exhq} hold with respect to
the entire (pruned) set $L$, and will make Lemma~\ref{le:line2flat}
hold for each 2-flat.

\noindent (d) Incidences between points and lines that are contained
in some irreducible component of $Z(f,g)$ that is a 2-flat. These
incidences will be analyzed explicitly below, using the properties
of flat points and lines, as presented in Section~\ref{se:algflat}.

This classification of incidences, especially those of types (c) and
(d), holds in general, for any polynomial $g$ satisfying the
properties assumed in this treatment (that it has degree $O(D)$ and
that $Z(f,g)$ is two-dimensional), and their treatment also applies
to these more general scenarios.
\paragraph{Incidences within hyperplanes and quadrics.}
We next derive a bound that we will use several times later on, in cases where
we can \emph{partition} $P$ and $L$ (or, more precisely, subsets thereof)
among some finite collection of hyperplanes
and quadrics, so that all the relevant incidences occur between points and
lines that are assigned to the same surface.
Recall that we have already applied a similar partitioning among the factors
of $f$ and of its derivatives.
The prime application of this bound will be to incidences of type (a) above,
but it will also be used in the analysis of type (d) incidences, and in the
analysis of the second case (b), where $Z(f,\flf)$ is three-dimensional, i.e.,
when $Z(f,\flf)=Z(f)$. In particular, we emphasize that the following
proposition does not require that $Z(f,\flf)$ be two-dimensional.
\begin {proposition}
\label {th:bwhyp} Let $H_1,\ldots,H_t$ be a finite collection of
hyperplanes and quadrics. Assume that the points of $P$ and the
lines of $L$ are partitioned among $H_1,\ldots, H_t$, so that each
point $p\in P$ (resp., each line $\ell\in L$) is assigned to a unique hyperplane
or quadric that contains $p$ (resp., fully contains $\ell$), and
assume further that each $H_i$ is $q$-restricted and that each 2-flat contains
at most $s$ lines of $L$. Then the overall number of incidences
between points and lines that are assigned to the same surface is
\begin {equation} \label{eq:hyps}
O\left(m^{1/2}n^{1/2}q^{1/4} + m^{2/3}n^{1/3}s^{1/3} + m + n\right).
\end {equation}
\end {proposition}

\noindent{\bf Proof.} For $i=1,\ldots,t$, let $L_{i}$ (resp.,
$P_{i}$) denote the set of lines of $L$ (resp., points of $P$), that
are assigned to $H_i$, and put $n_i := |L_{i}|$, $m_i := |P_{i}|$.
We have $\sum_i m_i = m$ and $\sum_i n_i = n$. For each $i$, since
$H_i$ is $q$-restricted, there exists a polynomial $g_i = g_{H_i}$,
of degree $O(\sqrt{q})$, such that all the lines of $L_i$, with the
exception of at most $q$ of them, are fully contained in ruled
components of $H_i\cap Z(g_i)$ that are not 2-flats. Write $L_i =
L_i^{nr}\cup L_i^r$, where $L_i^{r}$ is the subset of those lines
that are fully contained in ruled components of $H_i\cap Z(g_i)$
that are not 2-flats, and $L_i^{nr}$ is the complementary subset, of
size at most $q$. The lines in $L_i^r$ are contained in the union
$W^r$ of the ruled components of $H_i \cap Z(g_i)$ that are not
2-flats. We also remove from $L_i^r$ the subset $L_{i0}^r$ of $O(q)$
lines, as provided by Lemma~\ref{salta} (including all the lines
that are fully contained in more than one component $W_i$), and put
them in $L_i^{nr}$; we continue to use the same notations for these
modified sets. To apply Lemma~\ref{salta} to the case where $H_i$ is
a quadric, we first project the configuration onto some generic
3-space, and note that by Sharir and Solomon~\cite[Lemma
2.1]{surf-socg}, the projection of $W^r$ does not contain any
2-flat. Since the size of $L_i^{nr}$ is still $O(q)$, we have, by
Theorem~\ref{ttt},
\begin{align*}
I(P_i,L_i^{nr}) & = O\left( m_i^{1/2} |L_i^{nr}|^{3/4} +
m_i^{2/3} |L_i^{nr}|^{1/3}s^{1/3} + m_i + |L_i^{nr}| \right) \\
& = O\left( m_i^{1/2}n_i^{1/2}q^{1/4} +
m_i^{2/3}n_i^{1/3}s^{1/3} + m_i + n_i \right) .
\end{align*}
(Note that Theorem~\ref{ttt} is directly applicable when $H_i$ is a hyperplane,
and that it can also be applied when $H_i$ is a quadric, by projecting the
configuration onto some generic hyperplane, similar to what we have just noted
for the application of Lemma~\ref{salta}.)

We next bound $I(P_i,L_i^r)$, using Lemma~\ref{salta} (when $H_i$ is
a quadric, we apply it to the generic projection of $W^r$ to three
dimensions, as above). Since $\deg (g_i) = O(\sqrt{q})$, $W^r$ is of
degree $O(\sqrt q)$, and thus also its projection to three
dimensions (see, e.g.,~\cite{Har}), in case $H_i$ is a quadric. We
have already removed from $L_i^r$ the subset $L_{i0}^r$ provided by
the lemma, and so the lemma yields the bound
$$
I(P_i,L_i^r) = O\left(m_i^{1/2}n_i^{1/2}q^{1/4} + m_i + n_i\right) .
$$
%
%
That is, we have:
$$
I(P_i,L_i) = O\left( m_i^{1/2}n_i^{1/2}q^{1/4} +
m_i^{2/3}n_i^{1/3}s^{1/3} + m_i + n_i \right) .
$$
Summing these bounds for $i=1,\ldots,t$, and using H\"older's inequality (twice),
we get the bound asserted in (\ref{eq:hyps}). \proofend

\paragraph{The case where $f$ is linear or quadratic.}
(These are the cases $D=1,2$.) Let us apply
Proposition~\ref{th:bwhyp} right away to bound the number of
incidences when our (irreducible) $f$ is linear or quadratic, that
is, when $Z(f)$ is a hyperplane or a quadric.
Proposition~\ref{th:bwhyp} (together with assumption (i) of the
theorem) then implies the following bound.
\begin{equation} \label{ma:in12}
I(P,L) = O\left(m^{1/2}n^{1/2}q^{1/4} + m^{2/3}n^{1/3}s^{1/3} + m + n\right) ,
\end{equation}
which is subsumed by the main bound (\ref{ma:in}).


\paragraph{Incidences within 2-flats fully contained in $Z(f,g)$.}
\label{se:twoflats}
Assuming generic directions of the coordinate axes, we may assume that,
for every non-singular point $p\in P$, $T_p Z(f)$ is not orthogonal to any of the axes.
This allows us to use the flatness criterion developed in Section~\ref{se:algflat}
to each point of $P$.

As in previous steps of the analysis, we simplify the notation by
denoting the subsets of the points and lines that lie in the 2-flat
components of $Z(f,g)$ as $P$ and $L$, and their respective sizes as
$m$ and $n$. Each point $p\in P$ (resp., each line $\ell \in L$)
under consideration is contained (resp., fully contained) in at
least one of the 2-flats $\pi_1,\ldots,\pi_k$ that are fully
contained in $Z(f,g)$ (these are the linear irreducible components
of $Z(f,g)$, and we have $k=O(\Deg^2)$). Let $P^{(2)}$ (resp.,
$P^{(3)}$) denote the set of points $p\in P$ that lie in at most two
(resp., at least three) of these 2-flats. Assign each point $p \in
P^{(2)}$ to the (at most) two 2-flats containing it. Note that if $p
\in P^{(2)}$, then every line $\ell$ that is incident to $p$ can be
contained in at most two of the 2-flats $\pi_i$, and we assign
$\ell$ to those 2-flats. Let $L^{(2)}$ denote the set of lines
$\ell\in L$ such that $\ell$ is incident to at least one point in
$P^{(2)}$ (and is thus contained in at most two 2-flats $\pi_i$),
and put $L^{(3)} = L\setminus L^{(2)}$. For $i=1,\ldots,k$, let
$L_{i}^{(2)}$ (resp., $P_{i}^{(2)}$) denote the set of lines of
$L^{(2)}$ (resp., points of $P^{(2)}$), that are contained in
$\pi_i$, and put $n_i := |L_{i}^{(2)}|$, $m_i = |P_{i}^{(2)}|$.
(Note that we ignore here lines that are not fully contained in one
of these 2-flats; these lines are fully contained in other
components of $Z(f,g)$ and their contribution to the incidence count
has already been taken care of.) By construction,
$$
\sum_{i=1}^k m_i \le 2m, \quad \text{and}\quad \sum_{i=1}^k n_i \le 2n.
$$
Moreover, a point $p\in P^{(2)}$ can be incident only to lines of $L$ that are
contained in one of the (at most) two 2-flats that contain $p$, so we have
$I(P^{(2)}, L^{(2)}) \le \sum_{i=1}^k I(P_{i}^{(2)}, L_{i}^{(2)})$.
The Szemer\'edi--Trotter bound (\ref{inc2}) yields
\begin {equation}
\label {st:b} I(P_{i}^{(2)}, L_{i}^{(2)}) =
O\left(m_i^{2/3}n_i^{2/3}+m_i+n_i \right) , \quad\quad i =1,\ldots,k
.
\end {equation}
By assumption (ii) of the theorem, $n_i\le s$ for each
$i=1,\ldots,k$, so, summing over $i=1,\ldots,k$ and using H\"older's
inequality, we obtain
\begin {align}
\label {eq:P2L2}
I(P^{(2)}, L^{(2)}) \le \sum_{i=1}^k I(P_{i}^{(2)}, L_{i}^{(2)}) & =
O\left( \sum_{i=1}^k \Bigl( m_i^{2/3}n_i^{2/3} + m_i + n_i \Bigr) \right) \nonumber \\
& = O\left( \left(\sum_{i=1}^k m_i^{2/3}n_i^{1/3}s^{1/3}\right) + m + n \right) \\
& = O\left( \Bigl( \sum_{i=1}^k m_i \Bigr)^{2/3} \Bigl( \sum_{i=1}^k n_i \Bigr)^{1/3}s^{1/3} + m + n \right) \nonumber  \\
& = O\left( m^{2/3}n^{1/3}s^{1/3} + m + n \right) \nonumber  .
\end {align}
Consider next the points of $P^{(3)}$, each contained in at least
three 2-flats that are fully contained in $Z(f)$. All the points of
$P^{(3)}$ are linearly flat (see Section~\ref{se:algflat} for details),
and are therefore flat. Notice that each
such point can be incident to lines of $L^{(2)}$ and to lines
of $L^{(3)}$. We prune away each line $\ell \in L$ that
contains fewer than $3\Deg$ points of $P^{(3)}$, losing at most
$3n\Deg$ incidences in the process.

Each of the surviving lines contains at least $3\Deg-3$ flat points,
and is therefore flat, because the degrees of the nine polynomials
whose vanishing at $p$ captures the flatness of $p$, are all at most
$3\Deg-4$. In other words, we are left with the task of bounding the
number of incidences between flat points and flat lines. To simplify
this part of the presentation, we again rename the sets of these
points and lines as $P$ and $L$, and denote their sizes by $m$ and
$n$, respectively.

\paragraph{Incidences between flat points and lines.}
By Lemma \ref{le:lar}, all the (non-singular) points of a flat line
have the same tangent hyperplane. We assign each flat point $p \in
P$ (resp., flat line $\ell \in L$) to $T_p Z(f)$ (resp., to $T_p
Z(f)$ for some (any) non-singular point $p\in P\cap\ell$; again we
only consider lines incident to at least one such point). We have
therefore partitioned $P$ and $L$ among \emph{distinct} hyperplanes
$H_1,\ldots,H_t$, and we only need to count incidences between
points and lines assigned to the same hyperplane. By assumptions (i)
and (ii) of the theorem, the conditions of
Proposition~\ref{th:bwhyp} hold, implying that the number of these
incidences is
\begin{equation} \label{ma:fl}
O\left(m^{1/2}n^{1/2}q^{1/4} + m^{2/3}n^{1/3}s^{1/3} + m + n\right).
\end{equation}

As promised, after having bounded the number of incidences within
the 2-flats that are fully contained in $Z(f,g)$, we remove from $L$
the lines that are contained in such 2-flats, and continue the
analysis with the remaining subset.

\paragraph{In summary,}
combining the bounds in (\ref{eq:P2L2}) and (\ref{ma:fl}), Proposition~\ref{th:degsquared},
and Lemmas~\ref{le:exhq} and \ref{le:line2flat},
the overall outcome of the analysis for the first case is summarized in the
following proposition. (In the proposition, $f$ is one of the irreducible factors of
the original polynomial or of one of its derivatives, and $P$ and $L$ refer to the
subsets assigned to that factor.)
\begin{proposition} \label{pr:case1}
Let $g$ be any polynomial of degree $O(D)$ such that $Z(f,g)$ is two-dimensional,
let $P$ be a set of $m$ points contained in $Z(f,g)$, and let $L$ be a set of $n$
lines contained in $Z(f,g)$. Then
\begin{equation} \label{eq:case1}
I(P,L) = I(P^*,L^*) + O\left(
m^{1/2}n^{1/2}q^{1/4} + m^{2/3}n^{1/3}s^{1/3} + m + nD \right) ,
\end{equation}
where $P^*$ and $L^*$ are subsets of $P$ and $L$, respectively, so
that each hyperplane or quadric is $O(D^2)$-restricted with respect
to $L^*$, and each 2-flat contains at most $O(D)$ lines of $L^*$. We
also have the explicit estimate
\begin {equation}
\label {le:mines} I(P^*,L^*) = \min \{O\left(mD^2+nD\right), \
O\left(m+nD^4\right) \}.
\end {equation}
\end{proposition}

 \noindent{\bf Remarks.} (1) As already noted, lines that are
contained in 2-flats that are fully contained in $Z(f,g)$ have
already been taken care of, and thus do not belong to $L^*$, so the
application of Lemma~\ref{le:line2flat} shows that \emph{every}
2-flat contains only $O(D)$ lines of $L^*$, and the application of
Lemma~\ref{le:exhq} shows that \emph{every} hyperplane or quadric is
$O(D^2)$-restricted.

\noindent (2) When $m$ and $n$ are such that the bound on
$I(P^*,L^*)$ in (\ref{le:mines}) is dominated by
$O(m^{2/5}n^{4/5}+m)$, we use these bounds explicitly, and get the
induction-free refined bound in (\ref{ma:inx}). This remark will be
expanded and highlighted later, as we spell out the details of the
induction process.


\subsection*{Second case: $Z(f)$ is ruled by lines}


We next consider the case where the four-dimensional flecnode
polynomial $\fl{f}^4$ vanishes identically on $Z(f)$. By Theorem
\ref{th:flec2}, this implies that $Z(f)$ is ruled by (possibly
complex) lines.

In what follows we assume that $\Deg \ge 3$ (the cases $D=1,2$ have
already been treated earlier, using Proposition~\ref{th:bwhyp}).
We prune away points $p \in P$, with $|\Sigma_p|\le 6$ (the number of incidences
involving these points is at most $6m = O(m)$). For simplicity of
notation, we still denote the set of surviving points by $P$. Thus
we now have $|\Sigma_p| > 6$, for every $p\in P$.

Recalling the properties of the $u$-resultant of $f$ (that is, the
$u$-resultant associated with $F_1(p;v)$, $F_2(p;v)$, $F_3(p;v)$),
as reviewed in Section~\ref{se:algres}, we have, by
Corollary~\ref{coren}, that $U(p;u_0,u_1,u_2,u_3)\equiv 0$ (as a
polynomial in $u_0,\ldots, u_3$) for every $p \in P$.

We will use the following theorem of Landsberg, which generalizes
Theorem \ref{th:flec2}. It is stated here in a specialized and
slightly revised form, but still for an arbitrary hypersurface in
any dimension, and for any choice of the parameter $k$. Recall that
$\Sigma^k$ is the union of $\Sigma_p^k$ over all $p \in X$, namely,
it is the set of all lines that osculate to $Z(f)$ to order three at
some point on $Z(f)$.

The actual application of the theorem will be for $X=Z(f)$ (and
$d=4, k=3$). We refer the reader to Section~\ref{se:alglineon} for
notations and further details.


\begin{theorem}[Landsberg~\protect{\cite[Theorem 3.8.7]{IL}}]
\label{anlan} Let $X \subset \mathbb P^{d}(\cplx)$ be a
hypersurface, and let $k \ge 2$ be an integer, such that there is an
irreducible component $\Sigma_{0}^k \subset \Sigma^k$ satisfying,
for every point $p$ in a Zariski open set $\mathcal O \subset Z(f)$,
$\dim \Sigma_{0,p}^k
> d-k-1,$ where $\Sigma_{0,p}^k$ is the set of lines in $\Sigma_0^k$ incident to $p$.
Then, for each point $p\in \mathcal O$, all lines in
$\Sigma_{0,p}^k$ are contained in $X$.
\end{theorem}

To appreciate the theorem, we note that, informally, lines through a
fixed point $p$ have $d-1$ degrees of freedom, and the constraint
that such a line osculates to $X$ to order $k$ removes $k$ degrees
of freedom, leaving $d-k-1$ degrees. The theorem asserts that if the
dimension of this set of lines is larger, for most points on $X$,
then these lines are fully contained in $X$. Note also that this is
a ``local-to-global'' theorem---the large dimensionality condition
has to hold at every point of some Zariski open subset of $Z(f)$,
for the conclusion to hold.

If $U(p;u_0,u_1,u_2,u_3)$ does not vanish identically (as a
polynomial in $u_0,u_1,u_2,u_3$) at every point $p \in Z(f)$, then
at least one of its coefficients, call it $c_U$, does not vanish
identically on $Z(f)$. In this case, as $U$ vanishes identically at
every point of $P$ (as a polynomial in $u_0,u_1,u_2,u_3$), it follows that
$P$ is contained in the two-dimensional variety $Z(f,c_U)$.
Since $c_U$ has degree $O(\Deg)$ in $x,y,z,w$ (by Theorem \ref{ren}),
we can proceed exactly as we did in the case where $Z(f,\fl{f}^4)$ was 2-dimensional.
That is, we obtain the bound (\ref{eq:case1}) in Proposition~\ref{pr:case1}, namely,
\begin{equation} \label{eq:case1x}
I(P,L) = I(P^*,L^*) + O\left( m^{1/2}n^{1/2}q^{1/4} +
m^{2/3}n^{1/3}s^{1/3} + m + nD \right) ,
\end{equation}
where $P^*$ and $L^*$ are subsets of $P$ and $L$, respectively, so
that each hyperplane or quadric is $O(D^2)$-restricted with respect
to $L^*$, and each 2-flat contains at most $O(D)$ lines of $L^*$. We
also have the explicit estimate
$$
I(P^*,L^*) = \min \{O\left(mD^2+nD\right), \ O\left(m+nD^4\right)
\}.
$$

Therefore, since this case does not require the following analysis,
it suffices to consider the complementary situation, where we
assume that $U(p;u_0,u_1,u_2,u_3) \equiv 0$ at every point
$p \in Z(f)$ (as a polynomial in $u_0,u_1,u_2,u_3$). By Theorem
\ref{ren}, $\Sigma_p^3$ is infinite, so its dimension is positive,
for each such $p$.

Informally, the analysis proceeds as follows. Since $\Sigma_p^3$ is
(at least) one-dimensional for every point $p\in Z(f)$, the set
$\Sigma^3$, which is the union of $\Sigma_p^3$, over all $p\in
Z(f)$, has (at least) three degrees of freedom---three for
specifying $p$, at least one for specifying the line in
$\Sigma_p^3$, and one removed because the same line may arise at
each of its points (if it is fully contained in $Z(f)$). In what
follows we show that we can find a single irreducible component
$\Sigma_0^3$ of $\Sigma^3$, which is three-dimensional, and such
that for any point $p\in Z(f)$, the variety $\Sigma_{0,p}^3$ is at
least one-dimensional. This will facilitate the application of
Theorem~\ref{anlan} in our context.

\begin {theorem}
\label{th:cofa} There exists an irreducible component $\Sigma_0^3$
of $\Sigma^3$ of dimension at least three, such that for each
non-singular $p\in Z(f)$, the variety $\Sigma_{0,p}^3$ is at least
one-dimensional.
\end {theorem}

\noindent{\bf Proof.}
The proof makes use of the Theorem of the Fibers and related
results, as reviewed in Section~\ref{se:finfib}. Put
$$
W:=\{(p,\ell)\mid p \in \ell, \ell \in \Sigma_p^3\}\subset Z(f)
\times \Sigma^3.
$$
Note that $W$ is naturally embedded in $\P^3 \times \P^5$, where the
second component contains the Pl\"ucker hypersurface of lines in
3-space. $W$ can formally be defined as the zero set of homogeneous
polynomials; one polynomial defines the Pl\"ucker quadric, other
polynomials express the condition $p\in \ell$, and other polynomials
are those defining the projective variety $\Sigma_p^3$, whose
elements are now represented by their Pl\"ucker coordinates in the
appropriate projective space (see Section~\ref{se:alglineon} for
details). Therefore, $W$ is a projective variety.

Let $$\Psi_1 : W \to Z(f), \quad \Psi_2: W \to \Sigma^3$$ be the
(restrictions to $W$ of the) projections to the first and second
factors of the product. For an irreducible component $\Sigma_0^3$ of
$\Sigma^3$ (which is also a projective variety), put
$$
W_0:=\Psi_2^{-1}(\Sigma_0^3)=\{(p,\ell) \in W \mid \ell \in
\Sigma_{0,p}^3\}.
$$
Since $W$ and $\Sigma_0^3$ are projective
varieties, so is $W_0$. (Indeed, if $W=Z(\{f_i(p,\ell)\})$, and
$\Sigma_0^3=Z(\{g_j(\ell)\})$, for suitable sets of homogeneous
polynomials $\{f_i\}, \{g_j\}$, then
$W_0=Z(\{f_i(p,\ell),g_j(\ell)\})$.)

Let $\tilde W_0$ denote some irreducible component of $W_0$, and put
$Y:=\Psi_1(\tilde W_0)\subset Z(f)$. By the projective extension
theorem (see, e.g., Cox et al.~\cite[Theorem 8.6]{CLO}), $Y$ is also
a projective variety.

For a point $p\in Y$, the fiber of the map $\Psi_1|_{\tilde W_0}:
\tilde W_0 \to Y$ over $p$ is contained in $\{p\}\times
\Sigma_{0,p}^3=\{(p,\ell)\mid \ell \in \Sigma_{0,p}^3\}$ (this is
the fiber of $\Psi_1|_{W_0}$ over $p$, which clearly contains the
fiber of $\Psi_1|_{\tilde W_0}$ over $p$, as $\tilde W_0 \subseteq
W_0$).

We will show that there exists some component $\Sigma_0^3$, and some
irreducible component $\tilde W_0$ of $W_0=\Psi_2^{-1}(\Sigma_0^3)$,
such that (i) $Y=\Psi_1(\tilde W_0)$ is equal to $Z(f)$, and (ii)
for every point $p\in Z(f)$, the fiber of $\Psi_1|_{\tilde W_0}:
\tilde W_0 \to Y$ over $p$ is (at least) one-dimensional; in this
case we say that $\Sigma_0^3$ and $\tilde W_0$ form a
\emph{one-dimensional line cover} of $Z(f)$. Suppose that we have
found such a pair $\Sigma_0^3$, $\tilde W_0$. As noted above, the
fiber of $\Psi_1|_{\tilde W_0}$ over $p$ is contained in (or equal
to) $\{p\}\times \Sigma_{0,p}^3$, and $\dim(\{p\}\times
\Sigma_{0,p}^3)=\dim(\Sigma_{0,p}^3).$ Therefore, since $Y=Z(f)$,
this would imply that, for every $p\in Z(f)$, we have
$\dim(\Sigma_{0,p}^3)\ge 1$, which is what we want to prove.

We pick some component $\Sigma_0^3$, and some irreducible component
$\tilde W_0$ of $W_0=\Psi_2^{-1}(\Sigma_0^3)$, and analyze when do
$\Sigma_0^3$ and $\tilde W_0$ form a one-dimensional line cover of
$Z(f)$. Put, as above, $Y = \Psi_1(\tilde W_0)$.
For a point $p\in Y$, put $\lambda(p)=\dim(\Psi_1|_{\tilde
W_0}^{-1}(\{p\}))$, and let $\lambda = \min_{p \in Y} \lambda(p)$.
As noted above, $\lambda(p)\le \dim(\Sigma_{0,p}^3)$.

By the Theorem of the Fibers (Theorem~\ref{th:harr}), applied to the
map $\Psi_1|_{\tilde W_0}: \tilde W_0 \to Y \subseteq Z(f)$, we have
\begin {equation}
\label {eq:fiber} \dim(\tilde W_0)= \dim(Y)+\lambda.
\end {equation}

Observe that $\lambda \le 1$. Indeed, if $\lambda=2$, then there
exists some non-singular point $p\in Y$, such that $\Sigma_{0,p}^3$
is (at least) two-dimensional, implying that $Z(f)$ is a
three-dimensional cone; since $p$ is non-singular, $Z(f)$ is thus a
hyperplane, contrary to our assumptions.

Assume first that $Y=\Psi_1(\tilde W_0)$ is equal to $Z(f)$
(this is part (i) of the definition of a one-dimensional line cover).
If $\lambda = 1$, part (ii) of this property also holds,
and we are done. Assume then that $\lambda=0$.
By the first part of the Theorem of the Fibers (Theorem~\ref{th:harr}),
the subset $Y_1 = \{p\in Y \mid \lambda(p) \ge 1\}$ is Zariski closed
in $Y$, so it is a subvariety of $Z(f)$, of dimension
at most $2$. Hence, for
each $p$ in the Zariski open complement $Y\setminus Y_1$, the fiber
$\Psi_1|_{\tilde W_0}^{-1}(\{p\})$ is finite.

The remaining case is when $Y=\Psi_1(\tilde W_0)$ is properly contained
in $Z(f)$. Since $Z(f)$ is irreducible, $Y$ is of dimension at most two.

To recap, we have proved that for each component $\Sigma_0^3$ of $\Sigma^3$,
and each component $\tilde W_0$ of $W_0$, if the associated $Y$ is
properly contained in $Z(f)$, then the image of $\tilde W_0$ under
$\Psi_1$ (that is, $Y$) is at most two-dimensional; we refer to this
situation as being of the first kind. If $Y=Z(f)$ but $\lambda=0$ (these
are refered to as situations of the second kind), then, except for a two-dimensional
subvariety $Y_1$ of $Z(f)$, the fibers of the map $\Psi_1|_{\tilde W_0}$ are finite.

However, in the case under consideration, we have argued that, for
any non-singular point $p\in Z(f)$, the fiber
$\Psi_1^{-1}(p)=\{p\}\times \Sigma_p^3$ is (at least)
one-dimensional.

We apply this analysis to all the irreducible components $\Sigma_0^3$ of
$\Sigma^3$, and to all the irreducible components of the corresponding
$W_0=\Psi_2^{-1}(\Sigma_0^3)$. Let $Y^*$ denote the union of all the images
$Y$ of the first kind, and of all the excluded subvarieties $Y_1$ of the
second kind. Being a finite union of two-dimensional varieties,
$Y^*$ is two-dimensional.

The union, over the irreducible components $\Sigma_0^3$ of $\Sigma^3$,
of all the corresponding components $\tilde W_0$, covers $W$, and therefore,
for any non-singular point $p\in Z(f)$, the union over all the
components $\tilde W_0$ of the fibers of $\Psi_1|_{\tilde W_0}$ over
$p$ is equal to the fiber of $\Psi_1$ over $p$, which is
one-dimensional (and thus infinite).

We claim that there must exist some irreducible component $\Sigma_0^3$
of $\Sigma^3$, and a corresponding irreducible component $\tilde
W_0$ of $W_0$, such that $Y=\Psi_1(\tilde W_0)$ is equal to $Z(f)$,
and the corresponding $\lambda$ is equal to $1$. Indeed, if this were
not the case, take any non-singular point $p$ in $Z(f)\setminus Y^*$.
Since $p$ is not in the image $\Psi_1(\tilde W_0)$, for any $\tilde W_0$
of the first kind, the fiber of $\Psi_1|_{\tilde W_0}$ at $p$ is empty.
Similarly, since $p$ is not in the excluded set $Y_1$ for any $\tilde W_0$
of the second kind, the fiber of $\Psi_1|_{\tilde W_0}$ at $p$ is finite.
But then the fiber of $\Psi_1$ at $p$, being a finite union of (empty or)
finite sets, must be finite, a contradiction that establishes the claim.

Since for every $p \in Y=Z(f)$, $\lambda \le \lambda(p)\le
\dim(\Sigma_{0,p}^3)$, it follows that all the fibers
$\Sigma_{0,p}^3$ are (at least) one-dimensional, completing the
proof. \proofend

\noindent{\bf Remark.} One interesting corollary of the Theorem of
the Fibers is that if we know that for any point $p$ in a Zariski open
subset $\O$ of $Z(f)$, the fiber of $\Psi_1$ over $p$ (which is equal to
$\{p\}\times \Sigma_p^3$) is one-dimensional, then this is true
for the entire $Z(f)$. Indeed, by the Theorem of the Fibers
(Theorem~\ref{th:cofa}), the set
$\{p \in Z(f)\mid \dim(\Psi_1^{-1}(\{p\}))\ge 1\}$ is Zariski closed,
and, since it contains the Zariski open set $\O$, it must be equal to $Z(f)$.

By the preceding remark, Theorem \ref{anlan} (with $d=4, k=3$, $\O =
Z(f)$, and $\Sigma_0^3$ as specified by Theorem~\ref{th:cofa}) then
implies that $Z(f)$ is \emph{infinitely ruled} by lines, in the
sense defined in Section~\ref{se:algres}; that is, each point $p \in
Z(f)$ is incident to infinitely many lines that are fully contained
in $Z(f)$, and, moreover, $\Sigma_{0,p}^3 = \Sigma_{0,p}$ (which is
the set of lines in $\Sigma_0$ incident to $p$). That is, we have
shown that $\Sigma_0^3 = \Sigma_0$. In other words, for each $p\in
Z(f)$, $\Sigma_{0,p}$ is of dimension at least $1$, or,
equivalently, the cone $\Xi_{0,p}$ (which is the union of the lines
in $\Sigma_{0,p}$) is at least two-dimensional. If, for some
non-singular $p\in Z(f)$, the cone $\Xi_{0,p}$ were
three-dimensional, then, as already noted, $Z(f)$ would be a
hyperplane, contrary to assumption. Thus, for each non-singular
point $p\in Z(f)$, the cone $\Xi_{0,p}$ is two-dimensional, and
$\Sigma_{0,p}$ is one-dimensional. We also have
$\dim(\Sigma_0)=\dim(\Sigma_0^3)\ge 3$. We thus have

\begin {corollary}
\label {co:cofa2} The union of lines in $\Sigma_0^3=\Sigma_0$ is
equal to $Z(f)$, and $\dim(\Sigma_0)=\dim(\Sigma_{0}^3) \ge 3$.
\end {corollary}

\paragraph{Severi's theorem.}
The following theorem is a major ingredient in the present part of
our analysis. It has been obtained by Severi \cite{severi} in 1901,
and a variant of it is also attributed to Segre~\cite{Seg}; it is
mentioned in a recent work of Rogora \cite{Ro}, in another work of
Mezzetti and Portelli~\cite{MP}, and also appears in the unpublished
thesis of Richelson~\cite{Ric}. Severi's paper is not easily
accessible (and is written in Italian). As a small service to the
community, we sketch in Appendix A a proof of this theorem (or
rather of a special case of the theorem that arises in our context),
suggested to us by A.J.de Jong.

\begin {theorem}[Severi's Theorem~\protect{\cite{severi}}]
\label{th:sev} Let $X \subset \mathbb P^d(\cplx)$ be a
$k$-dimensional irreducible variety, and let $\Sigma_0$ be an irreducible
component of maximal dimension of $F(X)$, such that the lines of
$\Sigma_0$ cover $X$. Then the following holds.
\begin {enumerate}
    \item
If $\dim (\Sigma_0) = 2k - 2$, then $X$ is a copy of $\mathbb
P^k(\cplx)$ (that is, a complex projective $k$-flat).
    \item
If $\dim(\Sigma_0) = 2k -3$, then either $X$ is a quadric, or $X$ is
ruled by copies of $\mathbb P^{k-1}(\cplx)$, i.e., every point
\footnote{
  Similar to the definition in Section~\ref{se:algruled} for the case of lines, it suffices to require this property for every point in some Zariski-open subset of $X$.
  Here too one can show that the two definitions are equivalent. See also the companion paper~\cite[Lemma 11]{SS3dv}.} $p\in X$ is incident
to a copy of $\mathbb P^{k-1}(\cplx)$ that is fully contained in
$X$.
\end {enumerate}
\end {theorem}
As is easily checked, the maximum dimension of $\Sigma_0$ is $2k-2$.
Note also that the cases where $\dim \Sigma_0 < 2k-3$ are not
treated by the theorem (although they might occur); see
Rogora~\cite{Ro} for a (partial) treatment of these cases.

We apply Severi's theorem to $Z(f)$ and to the component $\Sigma_0$
obtained in Theorem~\ref{th:cofa} and Corollary~\ref{co:cofa2}, with
$k=3$ and with $\dim(\Sigma_0)=3=2k-3$. We thus conclude that either
$Z(f)$ is a quadric, a case ruled out in the present part of the
analysis (which assumed that $\deg(f) \ge 3$), or it is ruled by
2-flats.

\paragraph{The case where $Z(f)$ is ruled by $2$-flats.}
In the remaining case, every point $p \in Z(f)$ (see the footnote in
Theorem~\ref{th:sev}) is incident to at least one 2-flat
$\tau_p\subset Z(f)$. Let $D_p$ denote the set of 2-flats that pass
through $p$ and are contained in $Z(f)$.

For a non-singular point $p\in Z(f)$, if $\vert D_p \vert > 2$, then
$p$ is a (linearly flat and thus) flat point of $Z(f)$. Recall that
we have bounded the number of incidences involving flat points (and
lines) by partitioning them among a finite number of hyperplanes,
and by bounding the incidences within each hyperplane. (Recall that lines
incident to fewer than $3D-3$ points of $P$ have been pruned away, losing
only $O(nD)$ incidences, and that the remaining lines are all flat.)
Repeating this argument here, we obtain the bound
$$
O\left(m^{1/2}n^{1/2}q^{1/4} + m^{2/3}n^{1/3}s^{1/3} + m + n\right).
$$
In what follows we therefore assume that all points of $P$ are
non-singular and non-flat (call these points \emph{ordinary} for
short), and therefore $\vert D_p \vert = 1$ or $2$, for each such
$p$. Put $H_1(p)$ (resp., $H_1(p), H_2(p)$) for the 2-flat (resp.,
two 2-flats) in $D_p$, when $\vert D_p \vert = 1$ (resp., $\vert D_p
\vert = 2$).

Clearly, each line in $L$, containing at least one ordinary point $p
\in Z(f)$, is fully contained in at most two 2-flats fully contained
in $Z(f)$ (namely, the 2-flats of $D_p$).

Assign each ordinary point $p\in P$ to \emph{each} of the at most
two 2-flats in $D_p$, and assign each line $\ell \in L$ that is
incident to at least one ordinary point to the at most two 2-flats
that fully contain $\ell$ and are fully contained in $Z(f)$ (it is
possible that $\ell$ is not assigned to any $2$-flat---see below).
Changing the notation, enumerate these 2-flats, over all ordinary
points $p\in P$, as $U_1,\ldots,U_k$, and, for each $i=1,\ldots,k$,
let $P_i$ and $L_i$ denote the respective subsets of points and
lines assigned to $U_i$, and let $m_i$ and $n_i$ denote their
cardinalities. We then have $\sum_i m_i \le 2m$ and $\sum_i n_i \le
2n$, and the total number of incidences within the 2-flats $U_i$
(excluding lines not assigned to any 2-flat) is at most
$\sum_{i=1}^k I(P_i,L_i)$. This incidence count can be obtained
exactly as in the first case of the analysis, using the bound in
(\ref{eq:P2L2}). That is, we have
$$
\sum_{i=1}^k I(P_i,L_i) = O\left( m^{2/3}n^{1/3}s^{1/3} + m + n
\right) .
$$
As noted, this bound does not take into account incidences involving
lines which are not contained in any of the 2-flats $U_i$ (and are
therefore not assigned to any such 2-flat). It suffices to consider
only lines of this sort that are non-singular and non-flat, since
singular or flat lines are only incident to singular or flat points,
and we assumed above that all the points of $P$ are ordinary points.
If $\ell$ is a non-singular and non-flat line, and is not fully
contained in any of the $U_i$, we call it a \emph{piercing line} of
$Z(f)$.

\begin {lemma}
\label{le:pierc} If $\ell$ is a piercing line of $Z(f)$, then the
union of lines fully contained in $Z(f)$ and intersecting $\ell$ is
equal to $Z(f)$.
\end {lemma}

\noindent{\it Proof.} Let $V$ denote this union. By a suitable
extension to four dimensions of a similar result of Sharir and
Solomon~\cite[Lemma 5]{SS3d}, $V$ is a variety in the complex
projective setting, which we assume throughout this part of the
analysis. Clearly $V\subseteq Z(f)$. If $V$ is strictly contained in
$Z(f)$, then, since $Z(f)$ is irreducible, $V$ must be a finite
union of irreducible components $V_1,\ldots, V_k$, each of dimension
at most two. Let $p \in \ell$ be an ordinary point of $Z(f)$ (since
$\ell$ is non-singular and non-flat, such a point exists), and let
$H_1(p)$ be one of the at most two 2-flats in $D_p$. Note that
$H_1(p)$ is contained in $V$ (because it is a union of lines fully
contained in $Z(f)$ and intersecting $\ell$ at $p$). We claim that
there exists some $V_j$ such that $H_1(p)\subseteq V_j$. Indeed,
otherwise, the intersection $H_1(p)\cap V_j$ would be (at most)
one-dimensional for each $j=1,\ldots, k$ (a variety strictly
contained in a 2-flat is of dimension at most one), and therefore
$$
V\cap H_1(p)=\left(\bigcup_{j=1}^k V_j\right)\cap H_1(p) =
\bigcup_{j=1}^k \left(V_j\cap H_1(p)\right)
$$
is a finite union of varieties of dimension at most one,
contradicting the fact that $H_1(p)$ is contained in $V$ (and is of
dimension two). This contradiction establishes the claim. Since
$H_1(p)$ and $V_j$ are two-dimensional irreducible varieties and
$H_1(p)\subseteq V_j$, it follows that $H_1(p)=V_j$.

In other words, for each ordinary point $p\in \ell$ there exists a
2-flat $H_1(p)\in D_p$ which is equal to some component $V_j$.
Consider only the components $V_j$ that coincide with such a 2-flat.
Since there are only finitely many components $V_j$ of this kind,
one of them, call it $V_{j_0}$, has to intersect $\ell$ in
infinitely many points, and therefore $\ell \subseteq V_{j_0}$. That
is, $\ell$ is contained in the 2-flat $V_{j_0}$ that is fully
contained in $Z(f)$.

Now pick any ordinary point $p\in P\cap \ell$. By definition, since
$p\in V_{j_0}$, $V_{j_0}$ must be one of the (at most) two 2-flats
in $D_p$. But then $\ell$ is fully contained in that 2-flat, which
is one of the $U_i$'s, and therefore $\ell$ is not a piercing line.
This contradiction completes the proof. \proofend

\noindent {\bf Remark.} The last step of the proof shows that if a
non-singular and non-flat line $\ell$ contains a point of $P$ then
it is piercing (if and) only if it is not contained in any 2-flat
fully contained in $Z(f)$.

\begin {lemma}
\label{le:pierc2} Let $p\in Z(f)$ be an ordinary point. Then $p$ is
incident to at most one piercing line.
\end {lemma}

\noindent{\bf Proof.} Assume to the contrary that $p$ is incident to
two piercing lines $\ell_1,\ell_2 \in L$. We claim that the 2-flat
$\pi_{12}$ that is spanned by $\ell_1$ and $\ell_2$ is fully
contained in $Z(f)$ (and thus, by the preceding remark, $\ell_1$ and
$\ell_2$ are not piercing lines). Indeed, for any point $q \in
\ell_1$, Lemma~\ref{le:pierc} implies that there exists some line
$\ell_q \ne \ell_1$, incident to $q$, that intersect $\ell_2$ and is
fully contained in $Z(f)$. When $q$ varies along the non-singular
points of $\ell_1$, we get an infinite collection of lines, fully
contained in both $Z(f)$ and $\pi_{12}$, i.e., in their intersection
$Z(f)\cap \pi_{12}$. If $\pi_{12}$ is not contained in $Z(f)$ then
$Z(f)\cap \pi_{12} \ne \pi_{12}$ is a degree-$D$ plane curve, so by
Lemma~\ref{le:btl}, it contains at most $\Deg$ lines, and therefore
cannot contain the infinite union of lines $\bigcup_p \ell_p$.
\proofend

Therefore, each ordinary point $p\in P$ is incident to at most one
piercing line, and the total contribution of incidences involving
ordinary points and piercing lines is at most $m$.


\paragraph{In summary,}
combining the bounds that we have obtained for the various subcases
of the second case, we get the following proposition. As in the
first case, here $f$ refers to a single irreducible factor (of the
original polynomial or one of its derivatives), $D$ to its degree,
and $P$ and $L$ refer to the subsets of the original respective sets
of points and lines, that are assigned to $f$.
\begin{proposition} \label{pr:case2}
Let $P$ be a set of $m$ points contained in $Z(f)$, and let $L$ be a set of $n$
lines contained in $Z(f)$, and assume that $Z(f)$ is ruled by lines and that
$f$ is of degree $\ge 3$. Then
\begin{equation} \label{eq:case2}
I(P,L) = I(P^*,L^*) + O\left(
m^{1/2}n^{1/2}q^{1/4} + m^{2/3}n^{1/3}s^{1/3} + m + nD \right) ,
\end{equation}
where $P^*$ and $L^*$ are subsets of $P$ and $L$, respectively, so
that each hyperplane or quadric is $O(D^2)$-restricted with respect
to $L^*$, and each 2-flat contains at most $O(D)$ lines of $L^*$. We
also have the explicit estimate
\begin{equation} \label{plodd}
I(P^*,L^*) = \min \{O\left(mD^2+nD\right), \ O\left(m+nD^4\right)
\}.
\end{equation}
\end{proposition}

\subsection*{The induction} In summary, after having exhausted all possible
cases, we are in the following situation; we finally undo the
shorthand notations that we have used, and re-express the various
bounds in terms of the original parameters.

The first partitioning step has resulted in a collection of
irreducible polynomials, which we write as $f_1,\ldots,f_k$, with
respective degrees $D_1,\ldots,D_k$, all upper bounded by the degree
$D$ chosen in (\ref{le:mines}) for the original values of $m$ and
$n$. The points of $P$ have been \emph{partitioned} among the zero
sets $Z(f_1),\ldots,Z(f_k)$, into respective pairwise disjoint
subsets $P_1,\ldots,P_k$, including a leftover subset $P'$ of points
outside all the zero sets, and the lines of $L$ have been
\emph{partitioned} among the zero sets, into respective pairwise
disjoint subsets $L_1,\ldots,L_k$, so that the zero set to which a
line is assigned fully contains it, and including a leftover subset
$L'$ of lines not fully contained in any zero set. Put $m_i =
|P_i|$, $n_i = |L_i|$, for $i=1,\ldots,k$, and $m' = |P'|$, $n' =
|L'|$. Then $m' + \sum_{i=1}^k m_i = m$, and $n' + \sum_{i=1}^k n_i
= n$.

Then $I(P,L)$ is $I(P',L')+\sum_{i=1}^k I(P_i,L_i)$ plus the number
of incidences between points assigned to some $Z(f_i)$ and lines not
fully contained in $Z(f_i)$. (Note that $I(P\setminus P',L')$ also
counts incidences of this kind.) As we have argued, the total number
of these additional incidences is $O(nD)$. That is, we have, for
\emph{any} choice of the degree $D$,
\begin{equation} \label{ipart}
I(P,L) \le I(P',L') + O(nD) + \sum_{i=1}^k I(P_i,L_i) .
\end{equation}
For each $i$, the preceding analysis culminates in the following bound.
\begin{equation} \label{eq:all}
I(P_i,L_i) = I(P_{i}^*,L_{i}^*) + O\left(
m_i^{1/2}n_i^{1/2}q^{1/4} + m_i^{2/3}n_i^{1/3}s^{1/3} + m_i + n_iD \right) ,
\end{equation}
where, for each $i$, $P_{i}^*$ and $L_{i}^*$ are respective subsets
of $P_i$ and $L_i$, so that each hyperplane or quadric is
$O(D^2)$-restricted with respect to $L_{i}^*$, and each 2-flat
contains at most $O(D)$ lines of $L_{i}^*$. We also have the
explicit estimate
\begin{equation} \label{iplexp}
I(P_{i}^*,L_{i}^*) = \min \{O\left(m_iD^2+n_iD\right), \
O\left(m_i+n_iD^4\right) \}, \ \text{ for each }i.
\end{equation}
In addition, for the large values of $D$ in (\ref{eq:degofpart}), we
have
\begin{equation} \label{iptlt}
I(P',L') = O\left(m^{2/5}n^{4/5} + m \right) .
\end{equation}


\paragraph{Induction-free derivation of the bound.}
To proceed with the analysis, for general values of $m$ and $n$, we bound the various quantities
$I(P_{i}^*,L_{i}^*)$ using induction. However, as asserted in the theorems, the cases where
$m\le n^{6/7}$ or $m\ge n^{5/3}$ admit an induction-free argument that yields
the improved bound in (\ref{ma:inx}), and we first dispose of these cases.
(Recall that these are the original values of $m$ and $n$, the respective
sizes of the entire input sets $P$ and $L$.)

Assume first that $m\le n^{6/7}$. We substitute (\ref{eq:all}), the first bounds in
(\ref{iplexp}), and (\ref{iptlt}) into (\ref{ipart}). Using the
Cauchy-Schwarz and H\"older's inequalities, we have $\sum_i
m_i^{1/2}n_i^{1/2} \le m^{1/2}n^{1/2}$ and $\sum_i
m_i^{2/3}n_i^{1/3} \le m^{2/3} n^{1/3}$. We also have $\sum_i m_i
\le m$ and $\sum_i n_i \le n$. In total we thus get
\begin{align*}
I(P,L) & = O\left( m^{2/5}n^{4/5} + m + m^{1/2}n^{1/2}q^{1/4} +
m^{2/3}n^{1/3}s^{1/3} +
mD^2 + nD \right) \\
& = O\left( m^{2/5}n^{4/5} + m + m^{1/2}n^{1/2}q^{1/4} +
m^{2/3}n^{1/3}s^{1/3} + n\right) ,
\end{align*}
where we have used the fact that $mD^2+nD = O(m^{2/5}n^{4/5}+n)$ for
the choice $D=O(m^{2/5}/n^{1/5})$ in (\ref{le:mines}). This
establishes (\ref{ma:inx}) for this case. The case $m\ge n^{5/3}$ is
handled in the same manner, using the second bounds $O(m_i+n_iD^4)$ in (\ref{iplexp})
instead, and the fact that the sum of these bounds is $O(m)$ when
$m\ge n^{5/3}$.

\paragraph{The induction via a new partitioning.}
We now proceed with the general case, where induction is needed. To
simplify the notation, we (again, but only temporarily) drop the
indices, and consider one of many (possibly a nonconstant number of)
subproblems, involving a set $P$ ($=P_{i}^*$) of $m$ ($\le m_i$)
points and a set $L$ ($=L_{i}^*$) of $n$ ($\le n_i$) lines, so that
each hyperplane or quadric is $O(D^2)$-restricted for $L$, and each
2-flat contains at most $O(D)$ lines of $L$; here $D$ ($=D_i$) is
the degree of the corresponding factor $f$ ($=f_i$), which is upper
bounded by the value in (\ref{eq:degofpart}). In what follows we
will use this latter bound for (an upper bound on) the $D_i$'s.

To make the induction work, we choose a degree $E$, typically much
smaller than $D$ (see below for the actual value),
and construct a new partitioning polynomial $h$ of
degree $E$ for $P$. (Although $P \subset Z(f)$ and each
line of $L$ is fully contained in $Z(f)$, we ignore here $f$
completely, possibly losing some structural properties of $P$
and $L$, and consider only the partitioning induced by $h$.)
With an appropriate value of $r=\Theta(E^4)$, we obtain $O(r)$ cells,
each containing at most $m/r$ points of $P$, and each line of
$L$ either crosses at most $E+1$ cells, or is fully contained in $Z(h)$.

Set $P_0:= P\cap Z(h)$ and $P':=P\setminus P_0$. Similarly, denote
by $L_0$ the set of lines of $L$ that are fully contained in $Z(h)$,
and put $L':=L\setminus L_0$. We repeat the whole analysis done so
far, but with $h$ and its degree $E$ instead of $f$ and $D$, for the
points of $P$ and the lines of $L$. That is, we apply, to our $P$
and $L$, the bounds given in (\ref{ipart}), (\ref{eq:all}), and
(\ref{iplexp}) (but not the one in (\ref{iptlt})), with $E$ instead
of $D$. Moreover, in this application we exploit the property that
each hyperplane or quadric is $O(D^2)$-restricted with respect to
$L$, and each 2-flat contains at most $O(D)$ lines of $L$. We thus
get the following recurrence (where the parameters $k$, $P_i$,
$L_i$, etc., are new and depend on $h$, but we recycle the notation
in the interest of simplicity).
\begin{align*}
I(P,L) & \le I(P',L') + O(nE) + \sum_{i=1}^k I(P_i,L_i) \\
& = I(P',L') + O(nE) + \sum_{i=1}^k I(P_{i}^*,L_{i}^*) +
\sum_{i=1}^k O\left( m_i^{1/2}n_i^{1/2}D^{1/2} +
m_i^{2/3}n_i^{1/3}D^{1/3} + m_i + n_i E \right) .
\end{align*}

Concretely, $P'$ is the subset of the points of $P$ contained in the
cells of the $h$-partition, $L'$ is the subset of lines of $L$ not
fully contained in $Z(h)$, $P_i$ and $L_i$ are the subsets of the
points and lines assigned to the various irreducible factors $h_i$
of $h$ and of its derivatives, and $P_i^*, L_i^*$ are the excluded
subsets, as provided in Propositions~\ref{pr:case1}
and~\ref{pr:case2}.

Using the Cauchy-Schwarz and H\"older's inequalities in the second
sum, we get, for a suitable absolute constant $a$,
$$
I(P,L) \le I(P',L')
+ a \left( m^{1/2}n^{1/2}D^{1/2} + m^{2/3}n^{1/3}D^{1/3} + m + n E \right)
+ \sum_{i=1}^k I(P_{i}^*,L_{i}^*) .
$$

We have
$$
\sum_{i=1}^k I(P_{i}^*,L_{i}^*) \le a'\left( \sum_{i=1}^k \min\{
m_iE^2+n_iE, \ m_i+n_iE^4 \}\right) \le \min \{a'(mE^2+nE), \
a'(m+nE^4)\} ,
$$
for a suitable absolute constant $a'$.
That is, slightly increasing the coefficient $a$, we have
\begin{align} \label{eqind1}
I(P,L) \le I(P',L') + a \left( m^{1/2}n^{1/2}D^{1/2} +
m^{2/3}n^{1/3}D^{1/3} + m + n E \right) + \min\{amE^2, \ anE^4\}.
\end{align}

We next turn to bound $I(P',L')$. For each cell $\tau$ of
$\reals^4\setminus Z(h)$, put $P_\tau := P'\cap\tau$, and let
$L_\tau$ denote the set of the lines of $L'$ that cross $\tau$; put
$m_\tau = |P_\tau| \le m/r$ (where $r=\Theta(E^4)$), and $n_\tau =
|L_\tau|$. Since every line $\ell\in L'$ crosses at most $E+1$
components of $\reals^4\setminus Z(h)$, we have $\sum_\tau n_\tau
\le n(1+E)$.

To simplify the application of the induction hypothesis within the
cells of the partition, we want to make the subproblems be of
uniform size, so that $m_\tau = m/E^4$ and $n_\tau = n/E^3$ for each
$\tau$ (the latter quantity, up to some constant, is the average
number of lines crossing a cell). This is easy to enforce: To
achieve $m_\tau = m/E^4$, we simply partition $P_\tau$ into $\lceil
m_\tau/(m/E^4) \rceil = O(1)$ subsets, each consisting of at most
$m/E^4$ points, and analyze each subset separately. Similarly, if
$\tau$ is crossed by $\xi n/E^3$ lines, for $\xi>1$, we treat $\tau$
as if it occurs $\lceil \xi \rceil$ times, where each incarnation
involves all the points of (each of the constantly many
corresponding subsets of) $P_\tau$, and at most $n/E^3$ lines of
$L_\tau$. As is easily verified, the number of subproblems remains 
$O(E^4)$, with a larger constant of proportionality.


We apply the induction hypothesis for each cell $\tau$, to obtain
\begin{align*}
I(P_\tau,L_\tau) & \le 2^{c\sqrt{\log m_\tau}} \left(m_\tau^{2/5}n_\tau^{4/5} + m_\tau \right) +
\beta A\left( m_\tau^{1/2}n_\tau^{1/2}D^{1/2} + m_\tau^{2/3}n_\tau^{1/3}D^{1/3} + n_\tau \right) \\
& = 2^{c\sqrt{\log (m/E^4)}} \left( (m/E^4)^{2/5}(n/E^3)^{4/5} + m/E^4 \right) \\
& \quad + \beta A\left( (m/E^4)^{1/2}(n/E^3)^{1/2}D^{1/2} + (m/E^4)^{2/3}(n/E^3)^{1/3}D^{1/3} + n/E^3 \right) ,
\end{align*}
for a suitable absolute constant $\beta$. Summing this bound over
all cells $\tau$, that is, multiplying it by $O(E^4)$, we get, for a
suitable absolute constant $b$,
\begin{align} \label{cells1}
\sum_\tau I(P_\tau,L_\tau)
& \le b\cdot 2^{c\sqrt{\log (m/E^4)}} \left( m^{2/5}n^{4/5} + m \right) \\
& + b A\left( m^{1/2}n^{1/2}D^{1/2}E^{1/2} + m^{2/3}n^{1/3}D^{1/3}E^{1/3} + nE \right) \nonumber .
\end{align}
We have
\begin{align*}
2^{c\sqrt{\log (m/E^4)}} & =
2^{c\sqrt{\log m - 4\log E}} =
2^{c\sqrt{\log m} \left( 1 - \frac{4\log E}{\log m} \right)^{1/2} } \\
& < 2^{c\sqrt{\log m} \left( 1 - \frac{2\log E}{\log m} \right) }
= \frac{2^{c\sqrt{\log m}}} {2^{2c\log E / \sqrt{\log m}}}  .
\end{align*}
We choose $E$ to ensure that
$$
2^{2c\log E / \sqrt{\log m}} > 2 b , \quad\quad\text{or}\quad\quad
\frac{2c\log E }{ \sqrt{\log m}} > \log (2 b) ,
\quad\quad\text{or}\quad\quad \log E > \frac{\log (2 b)}{2c}
\sqrt{\log m} .
$$
That is, we choose
\begin{equation} \label{chE}
E > 2^{c^*\sqrt{\log m}} , \quad\quad\text{for}\quad\quad c^* =
\frac{\log (2 b)}{2c} < c/3 ,
\end{equation}
where the last constraint can be enforced if $c$ is chosen sufficiently large.
With this constraint on the choice of $E$, (\ref{cells1}) becomes
\begin{align} \label{cells2}
\sum_\tau I(P_\tau,L_\tau)
& \le \frac{1}{2} 2^{c\sqrt{\log m}} \left( m^{2/5}n^{4/5} + m \right) \\
& + b A\left( m^{1/2}n^{1/2}D^{1/2}E^{1/2} + m^{2/3}n^{1/3}D^{1/3}E^{1/3} + nE \right) \nonumber .
\end{align}
Adding this bound to the one in (\ref{eqind1}), we get
\begin{align} \label{allinc}
I(P,L)
& \le \frac{1}{2} 2^{c\sqrt{\log m}} \left( m^{2/5}n^{4/5} + m \right) \nonumber \\
& + (bA+a)\left( m^{1/2}n^{1/2}D^{1/2}E^{1/2} + m^{2/3}n^{1/3}D^{1/3}E^{1/3} + nE \right) + am \\
& +\min \{a m E^2, \ a n E^4\}. \nonumber
\end{align}

Returning to the original notations, we have just bounded
$I(P_{i}^*,L_{i}^*)$, for any $i=1,\ldots, k$. Concretely, we have
shown that, for each $i$,
\begin{align} \label{allinc2}
I(P_i^*,L_i^*)
& \le \frac{1}{2} 2^{c\sqrt{\log m_i}} \left( m_i^{2/5}n_i^{4/5} + m_i \right) \nonumber \\
& + (bA+a)\left( m_i^{1/2}n_i^{1/2}D^{1/2}E_i^{1/2} + m_i^{2/3}n_i^{1/3}D^{1/3}E_i^{1/3} + n_iE_i \right) + am_i \\
& +  \min \{a m_i E_i^2, \ a n_i E_i^4\},\nonumber
\end{align}
where $E_i$ is the degree of the new partitioning polynomial that is
constructed for $P_i^*$ and $L_i^*$.

We now add up these bounds, using~(\ref{ipart}), (\ref{eq:all}), and (\ref{iptlt}),
and replacing the $E_i$'s by a common upper bound $E$ that we will
choose shortly. We thus get the following bound, where now $P$ and $L$ stand,
respectively, for the original, entire input sets of points and lines.
\begin{align} \label{allinc3}
I(P,L)
& \le \gamma\left( m^{2/5}n^{4/5} + m \right) + \gamma nD  \nonumber \\
& + \gamma \sum_{i=1}^k \left( m_i^{1/2}n_i^{1/2} q^{1/4} +
m_i^{2/3}n_i^{1/3}s^{1/3} + m_i + n_i D \right) + \frac{1}{2}
\sum_{i=1}^k 2^{c\sqrt{\log m_i}} \left( m_i^{2/5}n_i^{4/5} + m_i
\right) \nonumber \\
& + \gamma \sum_{i=1}^k \left( m_i^{1/2}n_i^{1/2}D^{1/2}E^{1/2} + m_i^{2/3}n_i^{1/3}D^{1/3}E^{1/3} + n_iE + m_i \right) \\
& + \sum_{i=1}^k \min \{a m_i E^2, \ a n_i E^4\},\nonumber \\
\end {align}
for a suitable absolute constant $\gamma$.
With several applications of the Cauchy-Schwarz and H\"older's inequalities we get
\begin{align} \label{allinc4}
I(P,L)
& \le \left( \gamma+\frac 1 2 2^{c\sqrt {\log m}}\right)\left(
m^{2/5}n^{4/5} + m \right) \\
& + \gamma \left( m^{1/2}n^{1/2}q^{1/4} + m^{2/3}n^{1/3}s^{1/3} + m + n D \right) \nonumber \\
& + \gamma \left( m^{1/2}n^{1/2}D^{1/2}E^{1/2} +
m^{2/3}n^{1/3}D^{1/3}E^{1/3} + nE + m \right) + \min \{a m E^2, \ a
n E^4\} . \nonumber
\end{align}

We now bifurcate depending on the relation between $m$ and $n$,
where now, as in the recurrence just derived, $m$ and $n$ refer to
the original values of these parameters.

\paragraph{The case $m = O(n^{4/3})$.}
Recall that here we take $D=O(m^{2/5}/n^{1/5})$. It is easily
checked that, for this choice of $D$, each of the terms
$m^{1/2}n^{1/2}D^{1/2}$, $m^{2/3}n^{1/3}D^{1/3}$, $m$, and $nD\ge
n$, is $O(m^{2/5}n^{4/5})$, because $n^{1/2} \le m = O(n^{4/3})$.

We choose \footnote{This rather minuscule value of $E$ is only
needed when $m \approx n^{4/3}$; for smaller values of $m$, much
larger values of $E$ can be chosen.} $E = 2^{c^*\sqrt{\log m}}$.
This turns (\ref{allinc4}) into the bound
$$
I(P,L) \le \left(\gamma+\frac{1}{2} 2^{c\sqrt{\log m}} + \mu
2^{2c^*\sqrt{\log m}} \right) \left( m^{2/5}n^{4/5} + m \right) +
\gamma \left( m^{1/2}n^{1/2}q^{1/4} + m^{2/3}n^{1/3}s^{1/3}\right),
$$
for suitable absolute constants $\mu$ and $\gamma$. The choice of
$c^*$, and the assumption that $m\ge M_0$ and that $M_0$ is
sufficiently large, ensure that
$$
\gamma + \mu 2^{2c^*\sqrt{\log m}} < \frac{1}{2} 2^{c\sqrt{\log m}} ,
$$
and thus we get
$$
I(P,L) \le 2^{c\sqrt{\log m}} \left( m^{2/5}n^{4/5} + m \right)  +
\gamma \left( m^{1/2}n^{1/2}q^{1/4} + m^{2/3}n^{1/3}s^{1/3}\right),
$$
which is the bound asserted in (\ref{ma:in}).
\paragraph{The case $m = \Omega(n^{4/3})$.}
Here we take $D=O(n/m^{1/2})$. It is easily checked that, for this
choice of $D$, each of the terms $m^{1/2}n^{1/2}D^{1/2}$,
$m^{2/3}n^{1/3}D^{1/3}$, $m^{2/5}n^{4/5}$, and $nD\ge n$, is $O(m)$,
because $m = \Omega(n^{4/3})$.

We choose, as before, $E = 2^{c^*\sqrt{\log m}}$ (or a larger value when applicable), 
and note that, for $m\ge M_0$ sufficiently large, the term $nE^4$ is also $O(m)$. This
turns (\ref{allinc4}) into the bound
$$
I(P,L) \le \left(\gamma+\frac{1}{2} 2^{c\sqrt{\log m}} + \mu
2^{2c^*\sqrt{\log m}} \right) \left( m^{2/5}n^{4/5} + m \right) +  +
\gamma \left( m^{1/2}n^{1/2}q^{1/4} + m^{2/3}n^{1/3}s^{1/3}\right),
$$
for suitable absolute constants $\mu$ and $\gamma$. As above, the
choice of $c^*$, and the assumption that $m\ge M_0$ and that $M_0$
is sufficiently large, ensure that
$$
\gamma+ \mu 2^{2c^*\sqrt{\log m}} < \frac{1}{2} 2^{c\sqrt{\log m}} ,
$$
and thus we get
$$
I(P,L) \le 2^{c\sqrt{\log m}} \left( m^{2/5}n^{4/5} + m \right)  +
\gamma \left( m^{1/2}n^{1/2}q^{1/4} + m^{2/3}n^{1/3}s^{1/3}\right),
$$
again establishing the bound in (\ref{ma:in}). Therefore, in both
cases, we completed, at last, the induction step and thus
establishing the general upper bound (\ref{ma:in}) in the theorem.
The improved bound in (\ref{ma:inx}), for $m\le n^{6/7}$ or for
$m\ge n^{5/3}$, has already been established. With the lower bound
construction, given in the following section, the proof of the
theorem is completed. \proofend



\section{The lower bound} \label{sec:low}
In this section we present a construction that shows that the bound
asserted in the theorem is worst-case tight (except for the factor
$2^{c\sqrt{\log m}}$), for each $m$ and $n$, and for $q$ and $s$ in
suitable corresponding ranges, made precise below. 
The construction is a generalization
to four dimensions of a construction due to Elekes; see~\cite{El}.
(A three-dimensional generalization has been used in Guth and
Katz~\cite{GK2} for their lower bound construction.)

We have already remarked that the ``lower order'' terms
$m^{1/2}n^{1/2}q^{1/4}$ and $m^{2/3}n^{1/3}s^{1/3}$ are both
worst-case tight, as they can be attained by a suitable packing of
points and lines into hyperplanes (for the first term) or planes
(for the second term). Specifically, assume that $s\ge\sqrt{q}$, and
create $n/q$ parallel hyperplanes, and place on each of them $q$ 
lines and $mq/n$ points in a configuration that attains the 
three-dimensional lower bound as in Guth and Katz~\cite{GK2}.  
Note that in this construction no plane contains more than 
$\sqrt{q}\le s$ lines, as desired. Overall, we get
$$
(n/q)\cdot\Theta((mq/n)^{1/2}q^{3/4}) = \Theta(m^{1/2}n^{1/2}q^{1/4})
$$
incidences. A similar (and simpler) construction can be carried 
out for the second term $m^{2/3}n^{1/3}s^{1/3}$.

We therefore focus on the term $m^{2/5}n^{4/5}$ (the remaining terms
$m$ and $n$ are trivial to attain).

We fix two integer parameters $k$ and $\ell$, with concrete values
that will be set later, and take $P$ to be the set of vertices of
the integer grid
$$
\{ (x,y,z,w) \mid 1\le x\le k,\; 1\le y,z,w\le 2k\ell \} .
$$
We have $|P| = 8k^4\ell^3$.

We then take $L$ to be the set of all lines of the form
\begin{equation} \label{line}
y = ax+b, \quad\quad z = cx+d, \quad\quad w = ex+f ,
\end{equation}
where $1\le a,c,e\le \ell$ and $1\le b,d,f\le k\ell$. We have $|L| =
k^3\ell^6$. Note that each line in $L$ has $k$ incidences with the
points of $P$, one for each $x=1,2,\ldots,k$, so
$$
I(P,L) = k^4\ell^6 = \Theta(|P|^{2/5}|L|^{4/5}) ,
$$
as is easily checked. Note that $|L|^{1/2} \le |P| \le 8|L|^{4/3}$,
which is (asymptotically) the range of interest for this bound to be
significant: when $|P| < |L|^{1/2}$ we have the trivial bound
$I(P,L)=O(|L|)$, and when $|P|>|L|^{4/3}$, the leading term in the
bound changes qualitatively to $O(m)$, which is trivial for a lower
bound. Moreover, for any pair of integers $m$, $n$, with $n^{1/2}\le
m\le n^{4/3}$, we can find $k$ and $\ell$ for which $|P|=\Theta(m)$
and $|L|=\Theta(n)$. Specifically, choose $k=\Theta(m^{2/5}/n^{1/5})$
and $\ell = \Theta(n^{4/5}/m^{3/5})$; both are $\ge 1$ for the range
of $m$ and $n$ under consideration.

To complete the construction, we show that no hyperplane or quadric
can contain more than $q_0:= O\left( |L|^{6/5}/|P|^{2/5} \right) =
O(k^2\ell^6)$ lines of $L$, and no plane can contain more than
$s_0:=O\left( |L|^{7/5}/|P|^{4/5} \right)=O(kl^6)$ lines of $L$.
As an easy calculation shows, these threshold values of $q$ and $s$
are such that, for $q>q_0$ or $s>s_0$, the corresponding ``lower-dimensional''
term $m^{1/2}n^{1/2}q^{1/4}$ or $m^{2/3}n^{1/3}s^{1/3}$ dominates
the ``leading'' term $m^{2/5}n^{4/5}$ (for the former domination 
to arise, we need to assume, as above, that $\sqrt{q}\le s$), making the 
above construction pointless (see below for more details). The actual 
values of $q$ and $s$ that we will now derive are actually much smaller.

To estimate our $q$ and $s$, let $h$ be an arbitrary hyperplane.
If $h$ is orthogonal to the
$x$-axis then it does not contain any line of $L$, as is easily
checked, so we may assume that $h$ intersects any hyperplane of the
form $x=i$ in a $2$-plane $\pi_i$. The intersection of $P$ with
$x=i$ is a $2k\ell\times 2k\ell\times 2k\ell$ lattice, that we
denote as $Q_i$. Every line $\lambda\in L$ in $h$ meets $\pi_i$ at a
single point (as noted, it cannot be fully contained in $\pi_i$),
which is necessarily a point in $Q_i$ (every line of $L$ contains a
point of every $Q_i$). The size of $\pi_i\cap Q_i$ is easily seen to
be $O((k\ell)^2)$, and each point is incident to at most $\ell^2$
lines that lie in $h$. To see this latter property, substitute the
equations (\ref{line}) of a line of $L$ into the linear equation
defining $h$, say $Ax+By+Cz+Dw-1=0$ (where $B,C$ and $D$ are not
all 0). This yields a linear equation in $x$, whose $x$-coefficient
has to vanish. This in turn yields a linear equation in $a$, $c$,
and $e$, which can have at most $\ell^2$ solutions over
$[1,\ldots,\ell]^3$ (it is easily checked that the $x$-coefficient
cannot be identically zero for all choices of $a$, $c$, $e$). The
number of lines of point-line incidences of $P$ and $L$ within $h$
is thus $O(\ell^2(k\ell)^2)=O(k^2\ell^4)$. Since each line is
incident to $k$ points, necessarily all lying in $h$, it follows
that the number of lines of $L$ in $h$ is
$O(k^2\ell^4/k)=O(k\ell^4)$, which is always smaller than $q_0$.

This analysis easily extends to show that no quadric contains more
than $O(k\ell^4)$ lines of $L$; we omit the routine details.

Finally, let $\pi$ be a $2$-plane, where again we may assume that
$\pi$ is not orthogonal to the $x$-axis. Then $\pi$ meets a
hyperplane $x=i$ in a line $\mu$, and $\mu\cap Q_i$ contains at most
$k\ell$ points. Every line $\lambda$ in $\pi$ meets $\mu$ at one of
these points and, arguing as above, each such point can be incident
to at most $\ell$ lines that lie in $\pi$ (now instead of one linear
equation in $a$, $c$, $e$, we get two). Hence, $\pi$ contains at most
$k\ell^2/k=\ell^2$ lines of $L$, which is always smaller than $s_0$.

We have thus shown that the bound in Theorem~\ref{th:main} is (almost)
tight in the worst case. The bound will be tight when $|P|\le |L|^{6/7}$,
which occurs when $k\le \ell^{3/2}$, as an easy calculation shows.

\noindent{\bf Remark.}
As the analysis shows, the various constructions impose certain
constraints on the values of $q$ and $s$, and are therefore not as 
general (in terms of these parameters) as one might hope. It would 
be interesting to extend the constructions so that they apply to more 
general values of $q$ and $s$.

\section{Conclusion} \label{sec:conc}

The results of this paper (almost) settle the problem of point-line incidences
in four dimensions, but they raise several interesting and challenging open problems.
Among them are:

\smallskip

\noindent {\bf (a)} Get rid of the factor $2^{c\sqrt{\log m}}$ in
the bound. We have achieved this improvement when $m$ is not too
close to $n^{4/3}$, so to speak, allowing us to use the weak but
non-inductive bounds and complete the analysis in one step. We
believe that the ranges of $m$ where this can be done can be
enlarged, e.g. by improving the weak bounds. A concrete step in this
direction would be to improve the term $O(nD^4)$ in the second bound
in Proposition~\ref{th:degsquared}, which, as already remarked,
appears to be too weak. It would also be interesting to improve the
bound using the strategy in \cite{SSZ,SS3d}, which generates a
sequence of ranges of $m$, converging to $m=\Theta(n^{4/3})$, where
in each range the improved bound (\ref{ma:inx}) holds, with a
different constant of proportionality $A$. (For readers familiar with
the approaches in~\cite{SSZ, surf-socg}, we note that the reason
this technique does not appear to apply here is the multitude of
subproblems, each with its own $m_i, n_i$. The induction
in~\cite{SSZ, surf-socg} generates subproblems in which the relation
between $m$ and $n$ falls into a range already handled. Here though
we do not know how to enforce this property, as we have little
control over the values of $m,n$ in the resulting subproblems.

\smallskip

\noindent {\bf (b)} Extend (and sharpen) the bound of
Corollary~\ref{co:main} for any value of $k$. In particular, is it
true that the number of intersection points of the lines (this is
the case $k=2$; the intersection points are also known as
\emph{$2$-rich points}) is $O(n^{4/3}+nq^{1/2}+ns)$? We conjecture
that this is indeed the case. (In this conjecture we assume that we
have already managed to get rid of the factor $2^{c\sqrt{\log m}}$,
as in (a) above.)
A deeper question, extending a similar open problem in three
dimensions that has been posed by Guth and others (see, e.g., Katz's
expository note~\cite{Katz}), is whether the above conjectured bound
can be improved when $q=o(n^{2/3})$ and $s=o(n^{1/3})$, that is,
when the second and third terms in the conjectured bound become much
smaller than the term $n^{4/3}$. We also note that if we could
establish such a bound for the number of $k$-rich points, for any
constant $k$ (when $q$ and $s$ are not too large), then the case of
large $m$ (that is, $m=\Omega(n^{4/3})$) would become vacuous, as
only $O(n^{4/3})$ points could be incident to more than $k$ lines.

\smallskip

\noindent {\bf (c)} Extend the study to five and higher dimensions.
In a preliminary ongoing study, joint with Adam Sheffer, we can do
it using a constant-degree partitioning polynomial, with the
disadvantages discussed above (slightly weaker bounds, significantly
more restrictive assumptions, and inferior ``lower-dimensional''
terms). The leading terms in the resulting bounds, for points and
curves in $\reals^d$, are $O(m^{2/(d+1)+\eps}n^{d/(d+1)} +
m^{1+\eps})$, for any $\eps>0$. See also Dvir and Gopi~\cite{DG15} 
and Hablicsek and Scherr~\cite{HS14} for recent related studies.

Obtaining sharper results in such general settings, like the ones obtained in this paper, is
quite challenging algebraically, although some of the tools
developed in this work seem promising for higher dimensions too.

\smallskip

\noindent {\bf (d)} If we are given in advance that the points and
lines lie in some algebraic surface of a given degree $D>2$, can we
improve the bound and/or simplify the analysis? In our companion
work~\cite{SS3dv} we achieve these goals for the three-dimensional
case, improving the bound of Guth and Katz~\cite{GK2} in such
special cases.

\smallskip

\noindent
{\bf (e)} Elaborating on item (a) above, we note that the ``culprit''
Proposition~\ref{th:degsquared}, which produces the weak bounds that force us
to go into the induction, is only used in the case where $Z(f,g)$ is
two-dimensional, and the difficulty there lies in bounding the number of
incidences within a two-dimensional ruled surface (be it either one irreducible ruled surface
of large degree, or the union of many irreducible ruled surfaces of small degree). The analysis
of the three-dimensional analogous situation (addressed in Guth and Katz~\cite{GK2}),
cannot be applied here, since the degree of the underlying surface in four dimensions is $O(D^2)$
instead of $D$ in \cite{GK2}. In a recent study of Szermer\'edi-Trotter type
theorems in three dimensions~\cite{Kol}, Koll\'ar uses the \emph{arithmetic genus}
of curves to prove effective bounds on the number of point-line incidences
in three dimensions. In four dimensions, the situation is more involved,
but we hope that the arithmetic genus of the surface $Z(f,g)$ may yield
effective bounds for the number of incidences within this surface.


\paragraph{Acknowledgements.}
Work on this paper by Noam Solomon and Micha Sharir was supported by
Grant 892/13 from the Israel Science Foundation. Work by Micha
Sharir was also supported by Grant 2012/229 from the U.S.--Israel
Binational Science Foundation, by the Israeli Centers of Research
Excellence (I-CORE) program (Center No.~4/11), and by the Hermann
Minkowski-MINERVA Center for Geometry at Tel Aviv University. Part
of this research was performed while the authors were visiting the
Institute for Pure and Applied Mathematics (IPAM), which is
supported by the National Science Foundation. An earlier version of
this study appears in {\it Proc. 30th Annu. ACM Sympos. Comput.
Geom.}, 2014, 189--197, and the present version is also available in
arXiv:1411.0777v1.

We would like to thank several people whose advice, comments and
guidance have helped us a lot in our work on the paper. They are
J\'anos Koll\'ar, Martin Sombra, Aise J. de Jong, and Saugata Basu.
In addition, as noted, part of the work on the paper was carried out
during the special semester on Algebraic Techniques for
Combinatorial and Computational Geometry, held at the Institute for
Pure and Applied Mathematics at UCLA, in the Spring of 2014. We are
grateful for the pleasant working environment provided by IPAM, and
for the helpful interaction with additional colleagues, including
Larry Guth, Nets Hawk Katz, Terry Tao, Jordan Ellenberg, and many
others.

Author's address: School of Computer Science, Tel Aviv University,
Tel Aviv 69978, Israel; {\sl michas@tau.ac.il}, {\sl
noam.solom@gmail.com}.


\appendix

\section{Severi's Theorem}
\label{ap:sev}
In this appendix we sketch a proof of Severi's
theorem (Theorem \ref{th:sev}).

First, recall from Section~\ref{se:alglineon} that a real (resp.,
complex) surface $X$ is \emph{ruled by real (resp., complex) lines}
if every point $p\in X$ in a Zariski open dense set is incident to a
real (resp., complex) line that is fully contained in $X$. This
definition has been used in several recent works (see, e.g.,
~\cite{GK2}); this is a slightly weaker condition than the classical
condition that requires that \emph{every} point of $X$ be incident
to a line contained in $X$. Nevertheless, as we show next, the two
are equivalent.

\begin{lemma}
\label {le:rs} Let $f \in \reals[x,y,z] ($resp., $f \in
\reals[x,y,z,w])$ be an irreducible polynomial such that there
exists a Zariski open dense set $U\subseteq Z(f)$, so that each
point in the set is incident to a line, fully contained in $Z(f)$.
Then $\fl{f}\ ($resp., $\fl{f}^4)$ vanishes identically on $Z(f)$,
and $Z(f)$ is ruled by lines.
\end {lemma}
\noindent{\bf Proof.} By assumption and definition, $\fl{f}$ (resp.,
$\fl{f}^4$) vanishes on $U$. If it vanishes on $Z(f)$, Theorem
\ref{th:flec2} implies that $Z(f)$ is ruled. Otherwise, $Z(f,
\fl{f})$ (resp., $Z(f,\fl{f}^4)$) is properly contained in $Z(f)$
and contains $U$. Since $Z(f)$ is irreducible, this latter variety
must be of dimension at most $1$ (resp., $2$). On the other hand,
$Z(f, \fl{f})$ (resp., $Z(f,\fl{f}^4)$) is Zariski closed set (by
definition of the Zariski topology) and therefore contains its
Zariski closure. As $U$ is Zariski dense, its Zariski closure is
$Z(f)$. \proofend

\noindent{\bf Remark.} In Sharir and Solomon~\cite{SS3dv}, we have
proved the same statement without using the Flecnode polynomial.

This phenomenon generalizes to $k$-flats instead of lines (and the
proof translates verbatim).
\begin {lemma} \label{below}
Let $V$ be an irreducible variety for which there exists a Zariski
open subset $U\subseteq V$ with the property that each point $p\in
U$ is incident to a $k$-flat that is fully contained in $V$. Then
this property holds for every point of $V$.
\end {lemma}

We now proceed to sketch a proof of Severi's theorem. For convenience,
we repeat its statement.
\begin {thma} [Severi's Theorem~\protect{\cite{severi}}]
\label{th:sev2} Let $X \subset \mathbb P^d(\cplx)$ be a
$k$-dimensional irreducible variety, and let $\Sigma_0$ be a
component of maximal dimension of $F(X)$, such that the lines of
$\Sigma_0$ cover $X$. Then the following holds.
\begin {enumerate}
    \item
If $\dim (\Sigma_0) = 2k - 2$, then $X$ is a copy of $\mathbb
P^k(\cplx)$ (that is, a complex projective $k$-flat).
    \item
If $\dim(\Sigma_0) = 2k -3$, then either $X$ is a quadric, or $X$ is
ruled by copies of $\mathbb P^{k-1}(\cplx)$, i.e., every point $p\in
X$ is incident to a copy of $\mathbb P^{k-1}(\cplx)$ that is fully
contained in $X$.
\end {enumerate}
\end {thma}

We sketch a proof in the case $k=3$, $d=4$, under the simplifying
assumption that for any non-singular $x\in X$, $\Sigma_{0,x}$ is
infinite; this assumption holds in our application of the theorem
(by the informal dimensionality argument mentioned in the paper, it
holds ``on average'' in general for these parameters). Our proof is
based on a sketch provided by A.~J.~de Jong, via private
communication, and we are very grateful for his assistance.

\noindent{\bf Sketch of Proof.} For $x \in X$,  we recall that
$\Xi_{0,x}$ denotes the cone of lines (i.e., union of lines) of
$\Sigma_{0,x}$ The proof consists of the following steps.

\noindent (1) Assume first that $\dim (\Sigma_0)=2k-2=4$. Then there
exists some non-singular point $x_0\in X$ with
$\dim(\Sigma_{0,x_0})=2$. Indeed, if, for all non-singular points
$x\in X$, $\dim(\Sigma_{0,x}) \le 1$, then $\dim (\Sigma_0) < 4$
(see the analysis in Theorem~\ref{th:cofa}, and the preceding
analysis), contradicting the assumption in this case. By an argument
that has already been sketched earlier, this implies that
$\dim(\Xi_{0, x_0})=3$, i.e., the cone of lines in $\Sigma_{0,x_0}$
through $x_0$ is three-dimensional, and therefore $X=\Xi_{0,x_0}$.
As $x_0$ is non-singular, it follows that $X$ must be a hyperplane,
as claimed.

\noindent (2) Consider next the case where $\dim(\Sigma_0)=2k-3=3$,
and for any non-singular point $x\in X$, $\Sigma_{0,x}$ is
1-dimensional (as just argued, if $\Sigma_{0,x}$ is two-dimensional
for some non-singular $x \in X$, then $X$ is a hyperplane). In other
words, $\Sigma_{0,x}$, parameterized by the direction of its lines,
is a curve in $\mathbb P T_x X \cong \mathbb P^2(\cplx)$; put $e_x$
for its degree. If $e_x = 1$, then $\Xi_{0,x}$ contains a 2-flat.

We next define a ``plane-flecnode polynomial system'' associated
with $X$, that expresses, for a point $x\in X$, the existence of a
2-flat $H$, such that $H$ osculates to $X$ to order 3 at $x$. Since
$X$ is a hypersurface, we can write $X=Z(f)$, for a suitable
4-variate polynomial $f$ (see Section~\ref{se:alg}), and assume that
$f$ is irreducible (as $X$ is irreducible).

We represent a 2-flat through the origin in $\cplx^4$ (ignoring the
lower-dimensional family of 2-flats that cannot be represented in
this manner) as
\begin {equation} H_{v_0,v_1,v_2,v_3}:= \{(x,y,z,w) \mid  z=v_0x+v_1y, \  w=v_2x+v_3y\},
\end {equation}
 for $v_0,v_1,v_2,v_3 \in \cplx$.
The 2-flat $H_{v_0,v_1,v_2,v_3}$ is said to \emph{osculate to
$X=Z(f)$ to order $k$} at $p$, if the Taylor expansion of $f$ at $p$
along $H$ satisfies
\begin {equation}
f(p+(x,y,v_0x+v_1y, v_2x+v_3y))=O(x^{k+1}+y^{k+1}).
\end {equation}
This translates into a system of homogeneous polynomial equations in
$v_0,v_1,v_2,v_3$, involving the partial derivatives of $f$ up to
order $k$. Specializing to the case $k=3$, the plane-flecnode
polynomial system, $\pfl{f}$, associated with $f$, is obtained by
eliminating $v_0,v_1,v_2,v_3$ from these equations (for osculation
up to order $3$). This is the multipolynomial resultant system of
the polynomials defining these equations up to order $3$, with
respect to $v_0,v_1,v_2,v_3$ (see Van der Waerden~\cite[Chapter
XI]{Waer} for details).

Another theorem of Landsberg  \cite[Theorem 1]{Land} states that,
if, for every $g\in \pfl{f}$, $g$ vanishes identically on $X$, then
$X$ is ruled by 2-flats, which finishes the proof in this case.

Therefore, we may assume that $X\cap Z(\pfl{f})$ is a Zariski closed
proper subset of $X$. By definition of $\pfl{f}$, it follows that
for every non-singular point $x\in X \setminus Z(\pfl{f})$ (namely,
outside the Zariski closed set $Z(\pfl{f})$), we have $e_x>1$.
Indeed, if $e_x=1$, then, as observed above, there is a 2-flat
incident to $x$, and fully contained in $X$, implying that for every
$g \in \pfl{f}$, $g(x)=0$, contradicting the assumption that $x \in
X \setminus Z(\pfl{f})$.

For a generic hyperplane $H$ in $\mathbb P^4(\cplx)$, which is not
contained in $Z(\pfl{f})$, put $S_H := X \cap H$. As observed above,
$X \cap Z(\pfl{f})$ is properly contained in $X$, which in turn
implies that, for a generic hyperplane $H$ in $\mathbb P^4(\cplx)$,
$S_H$ is not fully contained in $Z(\pfl{f})$. Indeed, let $g$ be a
polynomial in $\pfl{f}$ that does not vanish identically on $X$.
Then $X \cap Z(g) = Z(f,g)$ is strictly contained in $X=Z(f)$, and
since $Z(f)$ is irreducible, it follows that $Z(f,g)$ is
two-dimensional. Therefore, for a generic hyperplane $H$, $X \cap
Z(\pfl{f}) \cap H$ is contained in the one-dimensional variety
$Z(f,g)\cap H$, and thus cannot contain the two-dimensional variety
$S_H$.

Let $x \in X$ be a non-singular point, and let $H$ be a hyperplane
in $\mathbb P^4(\cplx)$, which is incident to $X$ and not contained
in $Z(\pfl{f})$. We claim that for a generic $H$, there are $e_x$
\emph{distinct} lines that are incident to $x$ and fully contained
in $S_H$. Indeed, the intersection of the hyperplane $H$ with $T_x
X$ is a 2-flat in $T_x X$ containing $x$. Taking its
projectivization (where the point $x$ is regarded as $0$), namely,
$\mathbb P T_x X \cong \P^2$, the (generic) 2-flat $T_x X \cap H$
becomes a (generic) line. The degree of $\Sigma_{0,x} \subset
\mathbb P T_x X$ is $e_x$. Therefore, the intersection of $\Sigma_x$
with a line in $\mathbb P T_x X\cong \mathbb P^2$ consists of $e_x$
points, which are distinct since the line is generic. Therefore, its
intersection with $\Sigma_{0,x}$ consists of $e_x$ distinct points.
These $e_x$ distinct (projective) points represent $e_x$ distinct
lines, incident to $x$ and fully contained in $X\cap H=S_H$, as
claimed.

We say that a pair $(x,H)$, where $H$ is a hyperplane in
$\P^4(\cplx)$ and $x \in S_H$, is \emph{adequate} if there are $e_x$
distinct lines incident to $x$ that are fully contained in $S_H$.
Since a generic point $x$ is non-singular, the previous paragraph
implies that a generic pair $(x,H)$ is adequate. Therefore, by
changing the order of quantifiers, fixing a generic hyperplane $H$,
a generic point $x \in S_H$ is such that the pair $(x,H)$ is
adequate.

By Bertini's Theorem (see, e.g., Harris~\cite[Theorem 17.16]{Har}),
the irreducibility of $X$ implies that for a generic hyperplane $H$,
the surface $S_H$ is an irreducible surface in $H\cong \mathbb
P^3(\cplx)$. For a generic point $x \in S_H$, that is, outside an
algebraic curve $\mathcal C_H$ in $S_H$, the pair $(x,H)$ is
adequate. Therefore, there are $e_x$ distinct lines that are
incident to $x$ and fully contained in $S_H$, which, by Lemma
\ref{le:rs}, implies that $S_H$ is a ruled surface. Moreover, for
any $x \in S_H\setminus Z(\pfl{f})$, we have $e_x>1$. As observed
above, $\pfl{f}$ does not vanish identically on $S_H$, implying that
$Z(\pfl{f})\cap S_H$ is a Zariski closed proper subset of $S_H$,
i.e., an algebraic curve contained in $S_H$. Adding this curve to
$\mathcal C_H$, it follows that outside this algebraic curve, each
point of $S_H$ is incident to at least two lines fully contained in
$S_H$. By Sharir and Solomon~\cite[Lemma 9]{SS3dv}, this implies
that $S_H$ is either a 2-flat or a regulus. If $X$ is of degree
greater than two, then, for a generic hyperplane $H$, $S_H$ is a
(two-dimensional) surface of degree greater than two. Therefore, $X$
must be of degree at most two, namely, $X$ is either a hyperplane or
a quadric. If $X$ is a hyperplane, then $\Sigma_{0}$ is
four-dimensional, contrary to the present assumption, so finally, we
deduce that $X$ is a quadric, and the proof is complete. \proofend

\end{document}